\documentclass[12pt]{article}

\usepackage{latexsym,amsmath,amscd,amssymb,graphics}
\usepackage{enumerate}

\usepackage{graphicx}

\usepackage[square,authoryear]{natbib}
\usepackage[colorlinks]{hyperref}
\usepackage{url}

\usepackage[all]{xy}

\newcommand{\bfi}{\bfseries\itshape}

\makeatletter

\@addtoreset{figure}{section}
\def\thefigure{\thesection.\@arabic\c@figure}
\def\fps@figure{h, t}
\@addtoreset{table}{bsection}
\def\thetable{\thesection.\@arabic\c@table}
\def\fps@table{h, t}
\@addtoreset{equation}{section}

\makeatother

\begin{document}

\newtheorem{theorem}{Theorem}[section]
\newtheorem{definition}[theorem]{Definition}
\newtheorem{lemma}[theorem]{Lemma}
\newtheorem{remark}[theorem]{Remark}
\newtheorem{proposition}[theorem]{Proposition}
\newtheorem{corollary}[theorem]{Corollary}
\newtheorem{example}[theorem]{Example}

\def\below#1#2{\mathrel{\mathop{#1}\limits_{#2}}}



\title{The Lie-Poisson structure of the LAE-$\alpha$
equation}
\author{Fran\c{c}ois Gay-Balmaz$^{1}$ and Tudor
S. Ratiu$^{1}$}
\addtocounter{footnote}{1}
\footnotetext{Section de Math\'ematiques,
\'Ecole Polytechnique F\'ed\'erale de Lausanne.
CH--1015 Lausanne.
Switzerland.
\texttt{Francois.Gay-Balmaz@epfl.ch,
Tudor.Ratiu@epfl.ch}
\addtocounter{footnote}{1}
}

\date{  }

\maketitle

\makeatother

\maketitle




\begin{abstract}
This paper shows that the time $t$ map of the averaged Euler equations,
with Dirichlet, Neumann, and mixed boundary conditions is canonical
relative to a Lie-Poisson bracket constructed via a non-smooth
reduction for the corresponding diffeomorphism groups. It  is also shown
that the geodesic spray for Neumann and  mixed boundary conditions is
smooth, a result already known for Dirichlet boundary conditions.

\end{abstract}



\section{Introduction}
\label{section: Introduction}

The role of Hamiltonian structures for evolutionary
conservative equations in mathematical physics is well
established. In the finite dimensional case, that is, the
situation of ordinary differential Hamiltonian systems,
classical symplectic  and Poisson geometry and their
Lagrangian counterparts form the framework in which the
dynamics is formulated. When dealing with infinite
dimensional systems, namely the case of partial
differential equations, one is immediately confronted with
serious technical and conceptual difficulties. The main
issue is that, with the exception of certain
equations in quantum mechanics, all these PDEs need to be
formulated using a weak symplectic form. Also, for many
equations, the time evolution is not smooth in the
function spaces that are natural to the problem. If the
system is linear, this corresponds to the fact that the
right hand  side of the evolutionary equation is given by
an unbounded operator. Unfortunately, there is very little
general theory dealing  with the natural questions that
arise when working with Hamiltonian PDEs. The first
systematic attempt at such a devleopment can be found in
\cite{ChMa1974} and more recently, motivated by questions
regarding coherent states quantization, in
\cite{OdRa2003}. The present paper adds to this literature,
by presenting a precise Hamiltonian formulation of an
equation appearing in fluid dynamics.

\cite{Arnold1966} has given a Hamiltonian formulation of
the Euler equations for an incompressible homogeneous
perfect fluid (see also \cite{Arnold1989},
\cite{ArKh1998}, \cite{MaRa1999}). \cite{EbMa1970} have
shown that in appropriate Sobolev spaces, the Euler
equations are the spatial representation of a geodesic
spray that coincides with the dynamics of such a fluid in
material representation and that this geodesic spray is a
smooth vector field. In fact, this paper gives a rigorous
explanation with all the analytical details on how one
obtains the classical Euler equations as an
Euler-Poincar\'e equation associated to the group of
volume preserving  diffoemorphisms; the derivative loss
of the flow occuring in the passage from material to
spatial representation is also explained in this paper.
\cite{VaMa2004} have given a Hamiltonian formulation of
the Euler equations by carefully analyzing the function
spaces on which Poisson brackets are defined and
carrying out a Lie-Poisson reduction that takes into
account all analytical difficulties. They formulate an
analytical precise sense in which the flow of the Euler
equations are canonical. The remarkable fact is that the
passage from the previous analytically rigorous
Lagrangian formulation to this Hamiltonian picture is
nontrivial, mainly due to the fact that the flow is not
$C ^1$ from the Sobolev space of the initial condition to
itself. We shall comment below on the exact class of
Sobolev spaces needed in this formulation. A similar
analysis can be carried out for the incompressible
non-homogeneous Euler equations due to the resuls of
\cite{Marsden1976} which will involve semidirect product
groups.

The first goal of this paper is to carry out the program
outlined in \cite{VaMa2004}, that is, a non-smooth
Lie-Poisson reduction, for another equation appearing  in
fluid dynamics that has attracted a lot of attention
lately, namely the averaged or $\alpha$-Euler equation
(\cite{HoMaRa1998}). It has been shown in
\cite{MaRaSh2000}, \cite{Shkoller1998},
\cite{Shkoller2000} that these equations, either on
boundaryless manifolds or with Dirichlet
boundary conditions, have the same remarkable property,
namely in Lagrangian formulation they are smooth geodesic
sprays of $H ^1$-like weak Riemannian metrics on
appropriate diffeomorphism groups. These equations  are
intimately related to the Camassa-Holm equation
(\cite{CaHo1993}) for which this program can also be
carried out. We have chosen to work with the averaged
Euler equations because they have certain technical
difficulties not encountered for the homogeneous or
inhomogeneous Euler equations or the Camassa-Holm
equation; besides presenting more technical problems in
several steps, there also appears a one derivative loss
when formulating the precise sense in which  they are a
Lie-Poisson system and the flow is canonical.

The second goal of the paper is to show that the geodesic spray
for Neumann (or free-slip) and mixed boundary conditions is also smooth.
This completes the program outlined in \cite{MaRaSh2000}, \cite{Shkoller1998},
\cite{Shkoller2000} for these boundary conditions. This
shows in a different way that the averaged Euler
equations are well posed, a result due  to \cite{Shkoller2002} who uses
one  more derivative than  the present paper. We need this result in
order to achieve our third goal, namely to carry out a non-smooth
Lie-Poisson reduction for the averaged Euler equations with mixed
boundary conditions.

The plan of the paper is  as follows. Section 2 recalls
the relevant facts about the averaged Euler equations.
Section 3 gives the formulation of the averaged Euler
equations as a smooth geodesic spray of a weak Riemannian
metric on an appropirate group of volume preserving
diffeomeorphism. Section 4 gives the precise formulation
of the Poisson bracket, explicitly defines the correct
function spaces on which the Poisson bracket formula makes
sense and satisfies the usual axioms. Section 5 shows that
the averaged Euler equations are Hamiltonian relative to
the Poisson bracket defined previously with Hamiltonian
function given by the energy of the weak Riemannian
metric. It is also shown in what function spaces the flow
of these equations is a canonical map. The Lie-Poisson
reduction is also carried out explicitly in this section.
Section 6 proves the smoothness of the spray  for the averaged Euler
equations with mixed boundary conditions and generalizes to this case
all the results previously obtained in for Dirichlet boundary conditions.

\medskip

We close this introduction by presenting the geometric
setting of this paper and briefly recalling some of the
key facts about the Euler equations. Let
$(M,g)$ be a
$C^{\infty}$, compact, oriented, finite dimensional
Riemannian manifold of dimension at least two with
$C^{\infty}$ boundary
$\partial M$. The Riemannian volume form on $M $ is
denoted by $\mu$ and the induced volume form on $\partial
M$ by $\mu_{\partial}$. Let
$\nabla$ be the covariant derivative of the Levi-Civita
connection on $M$.

Let $N $ be another smooth boundaryless manifold. Recall
that if $s>\frac{1}{2}\dim M$ then a map $\psi: M \rightarrow N $
is of class $H ^s$ if its local representative in any pair
of charts is of class $H ^s$ as a map between open sets
of $\mathbb{R}^{\dim M} $ and $\mathbb{R}^{\dim N}$
respectively. If $s \leq \frac{1}{2}\dim M$ then, in general, a
map could be $H ^s$ in one pair of charts and fail to be
$H ^s$ in another one. Denote by $H^s(M, N) : = \{\psi: M
\rightarrow N \mid \psi \; \text{of class} \; H^s\}$ the
space of $H ^s$ maps from $M $ to $N $ for $s>\frac{1}{2}\dim M$. The set
$H^s(M,N)$ can be endowed with a smooth manifold structure
(see, e.g.,  \cite{EbMa1970, Palais1968}).

Let $\widetilde{M}$ denote the boundaryless double of
$M$. Then if $s> \frac{1}{2}\dim M + 1 $ the set
\begin{equation}
\label{diffeo_group}
\mathcal{D}^s:= \{ \eta \in H ^s(M ,\widetilde{M})
\mid \eta : M \rightarrow M \; \text{bijective},\;
\eta^{-1}
\in H ^s(M,
\widetilde{M}) \}
\end{equation}
is a group and a smooth submanifold of $H ^s(M,
\widetilde{M})$. If $\partial M = \varnothing$, then
$\mathcal{D}^s$ is an open subset of $H ^s(M, M)$. By the
Sobolev embedding  theorem, $\eta \in \mathcal{D}^s$ and
its inverse are necessarily of class $C ^1$. Therefore,
$\eta(\partial M) \subset \partial M $. The tangent space at the identity
$T_e\mathcal{D}^s$ consists of the $H^s$ class vector fields on $M$ which are
tangent to $\partial M$, denoted by
$\mathfrak{X}^s_{|\!|}$. Let
\begin{equation}
\label{}
\mathcal{D}^s_{\mu}:=\{\eta\in\mathcal{D}^s\mid\eta^*\mu=\mu\}
\end{equation}
be the subset of $\mathcal{D}^s$ whose elements preserve $\mu$. As proven in
\cite{EbMa1970}, the set $\mathcal{D}^s_\mu$ is a subgroup and a smooth
submanifold of $\mathcal{D}^s$. The tangent space
$T_e\mathcal{D}^s_{\mu}$
at the identity  equals $\mathfrak{X}^s_{\operatorname{div},
 |\!|}: = \{u \in \mathfrak{X}_{|\!|}^s
\mid \operatorname{div} u = 0\}$, the vector space of all
$H^s$ divergence free vector fields tangent to the
boundary. If $\dim M = 1 $ each of its connected components
is diffeomorphic to the circle $S ^1$. Taking on $S ^1$ the
usual length function, we see that the volume preserving
diffeomorphisms on the  circle are rotations. So, in this
case we have for each connected component
$\mathcal{D}^s_{\mu} = S ^1 $, which is not  an
interesting case.   Thus, since $\dim M
\geq 2$ we always have $s > 2 $.

On $\mathfrak{X}^s$ we can introduce the $L ^2 $ inner
product
\[
\langle u,v\rangle_0:=\int_Mg(x)(u(x),v(x))\mu(x)
\]
for any $u, v \in \mathfrak{X}^s$. This inner product on
$\mathfrak{X}^s$ is the value at the identity of two
distinct weak Riemannian metrics on $\mathcal{D}^s$, namely
\[
\mathcal{G}^0(\eta)(u_\eta,v_\eta)
:=\langle u_\eta\circ
\eta^{- 1},v_\eta\circ \eta^{-1}\rangle_0
\]
and
\[
\mathcal{G}(\eta)(u_\eta,v_\eta)
:= \int_M g(\eta(x))(u_\eta(x), v_\eta(x)) \mu(x)
\]
for any $u_\eta, v_\eta \in T_\eta \mathcal{D}^s$. Note
that $\mathcal{G}^0 $ is right invariant by construction,
whereas $\mathcal{G}$ is not. Their pull backs to
$\mathcal{D}^s_\mu$ coincide and yield a right invariant
weak Riemannian metric on $\mathcal{D}^s_\mu$. The
Euler equations
\begin{align*}
\label{Euler equations}
&\partial_t u(t) + \nabla_{u(t)} u(t) = -
\operatorname{grad} p(t)\\
& u(t) \in \mathfrak{X}^s_{\operatorname{div}, |\!|} \,,
\qquad u(0) = u_0
\quad \text{given}
\end{align*}
are the spatial representation of the geodesic spray on
$\mathcal{D}^s_\mu$ relative to this weak Riemannian metric
on $\mathcal{D}^s_\mu$ and this geodesic spray is a
smooth vector field on $T\mathcal{D}^s_\mu$ (see
\cite{EbMa1970}). The averaged Euler equations will be
presented in the next section.

\section{The geometry of LAE-$\alpha$ equation}

In this section we shall quickly review the results of
\cite{Shkoller2000} regarding the motion of the averaged
Euler equations. For $s>1 + \frac{1}{2}\dim
M$ we define three subsets of $\mathcal{D}^s$ which correspond to various boundary conditions.
The {\bfi Dirichlet diffeomorphism group\/} is defined by
\[
\mathcal{D}^s_D:=\{\eta\in\mathcal{D}^s\mid
\eta_{|\partial M}=id_{\partial M}\}.
\]
The {\bfi Neumann diffeomorphism group\/} is defined by
\[
\mathcal{D}^s_N:=\{\eta\in\mathcal{D}^s\mid
(T\eta_{|\partial M}\circ n)^{tan}=0 \text{ on $\partial M$}\},
\]
where $n$ denotes the outward-pointing unit normal
vector field along the boundary $\partial M$, and $(\cdot)^{tan}$ denotes the tangential part to the  boundary of a vector in $TM|\partial M$. The {\bfi mixed diffeomorphism group\/} is defined by
\[\mathcal{D}^s_{mix}:=\{\eta\in\mathcal{D}^s\mid \eta \text{ leaves $\Gamma_i$ invariant},\,\eta_{|\Gamma_1}=id_{|\Gamma_1},\,(T\eta_{|\Gamma_2}\circ n)^{tan}=0 \text{ on $\Gamma_2$}\},\]
where $\Gamma_1$ and $\Gamma_2$ are two disjoint subsets of $\partial M$ such that $\partial M=\Gamma_1\cup\Gamma_2$ and $\overline{\Gamma}_1=\partial M\setminus\Gamma_2$; furthermore, we assume that for all $m\in \Gamma_i$ we can find a local chart $U$ of $M$ at $m$ such that $\overline{U}\cap \partial M\subset\Gamma_i$.

The groups $\mathcal{D}^s_D, \mathcal{D}^s_N$ and $\mathcal{D}^s_{mix}$ are smooth Hilbert submanifolds and subgroups
of $\mathcal{D}^s$. The corresponding tangent spaces at the identity are given by
\[
\mathcal{V}^s_D
:=T_{id_M}\mathcal{D}^s_D=\{u\in\mathfrak{X}^s_{|\!|} \mid \,u_{|\partial M}=0\},
\]
\[
\mathcal{V}^s_N
:=T_{id_M}\mathcal{D}^s_N=\{u\in\mathfrak{X}^s_{|\!|} \mid (\nabla_nu_{|\partial M})^{tan}+S_n(u)=0 \text{ on $\partial M$}\},
\]
\[
\mathcal{V}^s_{mix}
:=T_{id_M}\mathcal{D}^s_{mix}=\{u\in\mathfrak{X}^s_{|\!|} \mid (\nabla_nu_{|\Gamma_1})^{tan}+S_n(u)=0 \text{ on $\Gamma_1$},\,u_{|\Gamma_2}=0\},
\]
where $S_n : T\partial M\rightarrow T\partial M$ is the Weingarten map
defined by $S_n(u):=-\nabla_un$. We can also form the corresponding sets $\mathcal{D}^s_{\mu,D}, \mathcal{D}^s_{\mu,N}$ and $\mathcal{D}^s_{\mu,mix}$ which have the volume-preserving constraint imposed. These sets are smooth Hilbert
submanifolds and subgroups of $\mathcal{D}^s_\mu$ and
$\mathcal{D}^s$. The corresponding tangent spaces at the identity are given by
\[
\mathcal{V}^s_{\mu,D}
:=T_{id_M}\mathcal{D}^s_{\mu,D}=\{u\in\mathfrak{X}^s_{\operatorname{div},|\!|} \mid \,u_{|\partial M}=0\},
\]
\[
\mathcal{V}^s_{\mu,N}
:=T_{id_M}\mathcal{D}^s_{\mu,N}=\{u\in\mathfrak{X}^s_{\operatorname{div},|\!|} \mid (\nabla_nu_{|\partial M})^{tan}+S_n(u)=0 \text{ on $\partial M$}\},
\]
\[
\mathcal{V}^s_{\mu,mix}
:=T_{id_M}\mathcal{D}^s_{\mu,mix}=\{u\in\mathfrak{X}^s_{\operatorname{div},|\!|} \mid (\nabla_nu_{|\Gamma_1})^{tan}+S_n(u)=0 \text{ on $\Gamma_1$},\,u_{|\Gamma_2}=0\}.
\]
Note that, as vector spaces,
$\mathcal{V}^r_D$ and $\mathcal{V}^r_{\mu,D} $ make sense
for $r \geq 1$, and $\mathcal{V}^r_N, \mathcal{V}^r_{mix}, \mathcal{V}^r_{\mu,N}$ and $\mathcal{V}^r_{\mu,mix}$ make sense for $r\geq 2$
but it is only for $s > 1+\frac{1}{2}\dim M$ that
they are the tangent spaces at the identity to the corresponding diffeomorphism subgroups.
If $1\leq r <2$ we set
\[
\mathcal{V}^r_N:=\mathfrak{X}^r_{|\!|},\qquad\mathcal{V}^r_{mix}
:=\{u\in\mathfrak{X}^r_{|\!|} \mid \,u_{|\Gamma_2}=0\}
\]
\[ \mathcal{V}^r_{\mu,N}:=\mathfrak{X}^r_{\operatorname{div},|\!|}, \qquad\mathcal{V}^r_{\mu,mix}
:=\{u\in\mathfrak{X}^r_{\operatorname{div},|\!|} \mid \,u_{|\Gamma_2}=0\}.
\]
For an arbitrary constant $\alpha> 0 $, consider on
$\mathfrak{X}^1$ the inner product
\begin{equation}
\label{h one inner product}
\langle
u,v\rangle_1:=\int_M\big{(}g(x)(u(x),v(x))+2
\alpha^2\overline{g}(x)(\operatorname{Def}
(u)(x),\operatorname{Def}(v)(x))\big{)}\mu(x),
\end{equation}
for all $u,v\in\mathfrak{X}^1$, where
\begin{equation}
\label{def operator}
\operatorname{Def}(u):=
\frac{\nabla u+(\nabla u)^t}{2}
\end{equation}
is the {\bfi deformation tensor\/}. In this formula,
$(\nabla u)^t$ denotes the transpose of the $(1,1)
$-tensor $\nabla u$ relative  to the metric $g $, that is,
$g(\nabla_v u, w) = g(v, (\nabla u)^t(w))$, for all $u, v,
w \in \mathfrak{X}^1 $. The symbol $\overline{g}$ denotes
the naturally induced inner product on $(1,1)$-tensors; in
coordinates, if $R, S $ are $(1,1)$-tensors then
$\overline{g}(R,S) = g_{ik} g^{j\ell} R^i_j S^k_\ell =\operatorname{Tr}
(R^t\cdot S)$.
This inner product induces by right translations a right
invariant weak Riemannian metric on
$\mathcal{D}^s_{\mu,mix}$ given by
\begin{equation}
\label{h one metric}
\mathcal{G}^1(\eta)(u_\eta,v_\eta):=\langle
u_\eta\circ
\eta^{- 1},v_\eta\circ \eta^{-1}\rangle_1
\end{equation}
for $u_\eta,v_\eta\in T_\eta\mathcal{D}^s_{\mu,mix}$.

We shall use throughout the paper the index
lowering and raising operators $\flat: \mathfrak{X}
\rightarrow \Omega ^1 $ and $\sharp : = \flat ^{-1}:
\Omega^1 \rightarrow  \mathfrak{X}$ induced by the metric
$g $, that is, $u ^\flat: = g (u, \cdot) $ for any $u \in
\mathfrak{X}$. Our conventions for the curvature and the
Ricci tensor and operator are
\[
\operatorname{R}(u,v): = \nabla_u \nabla_v -  \nabla_v \nabla_u  -
\nabla_{[u, v]}
\]\[
\overline{\operatorname{R}}(u,v,w,z):=g(\operatorname{R}(u,v)w,z)
\]\[
\operatorname{Ricci}(u,v):=\operatorname{Tr}(w\longmapsto\operatorname{R}(w,u)
v)
\]\[
g(\operatorname{Ric}(u), v) : = \operatorname{Ricci}(u,v)
\]

Let $\delta$ be the codifferential
associated to $g $. We denote by
\[
\Delta u=-[(d\delta+\delta d)u^{\flat}]^{\sharp}
\]
the usual {\bfi Hodge Laplacian\/} on vector fields and let
\[
\Delta_r:=\Delta+2\operatorname{Ric}
\]
be the {\bfi Ricci Laplacian}. We shall also need the operator
\[
\mathcal{L}:=\Delta_r+\operatorname{grad}\operatorname{div}.
\]
wich appears in the following formula (\cite {Shkoller2000})
\begin{equation}
\label{useful formula}
\langle u,v\rangle_1=\langle (1-\alpha^2\mathcal{L})u,v\rangle_0\text{ for
all }u,v\in\mathcal{V}^2_{mix}
\end{equation}
that will be used many times in this paper. For completeness we
shall provide below a  complete proof. Denote by
$\mathfrak{X}^{C^2}(U)$  the $C^2$ vector fields on an
open subset $U$ of $M$. We begin with the following.

\begin{lemma}$($Weitzenb\"ock formula$)$
\label{Weitzenbock}
Let $\{e_i \mid i = 1, \dots n\}$ be a local orthonormal
frame on an open subset
$U$ of $M$. Then on $\mathfrak{X}^{C^2}(U)$ the following
identity holds:
\begin{equation}
\label{Weitzenbock formula}
\Delta=\nabla^2_{e_i,e_i}-\operatorname{Ric}
\end{equation}
where $\nabla^2_{e_i,e_i}:=\nabla_{e_i}\nabla_{e_i}-\nabla_{\nabla_{e_i}e_i}$
is the second covariant derivative. In particular we
remark that $\nabla^2_{e_i,e_i}$ does not depend on the
local  orthonormal frame and so can be defined globally on
$M$.
\end{lemma}
\textit{Proof} : We will use the formula $\delta\alpha=-\operatorname{\bf i}_
{e_i}(\nabla_{e_i}\alpha)$ where $\{e_i\}$ is a local
orthonormal frame on an  open subset $U$ of $M$ and
$\alpha$ is a $k$-form (see \cite{Pe1997}). We also  need
the identities $d\alpha
(u,v)=(\nabla_u\alpha)(v)-(\nabla_v\alpha)(u)$  where
$\alpha$ is a one-form and $\nabla _u v^\flat = (\nabla _u
v)^\flat$ for any vector fields $u, v $ on $M $. Let $u \in
\mathfrak{X}^{C^2}(U)$ and recall that $\delta u^{\flat}=-
\operatorname{div}(u)$
. On $U$ we have :\\
\[
d(\delta u^{\flat})(v)=-d(\operatorname{div}(u))(v)=-d(g(\nabla_{e_i}u,e_i))(v)
=-g(\nabla_v\nabla_{e_i}u,e_i)-g(\nabla_{e_i}u,\nabla_ve_i).
\]
We also have:
\begin{eqnarray*}
\delta(d u^{\flat})(v)&=&-\operatorname{\bf i}_{e_i}(\nabla_{e_i}(d
u^{\flat}))
(v)=-\nabla_{e_i}(d u^{\flat})(e_i,v)\\
&=&-\nabla_{e_i}(d u^{\flat}(e_i,v))+d u^{\flat}(\nabla_{e_i}e_i,v)+d
u^{\flat}
(e_i,\nabla_{e_i}v)\\
&=&-\nabla_{e_i}\left(\nabla_{e_i}u^{\flat}(v)
- \nabla_vu^{\flat}(e_i)\right)+\nabla_
{\nabla_{e_i}e_i}u^{\flat}(v)-\nabla_vu^\flat(\nabla_{e_i}e_i)\\
&&+\nabla_{e_i}u^{\flat}(\nabla_{e_i}v)-\nabla_{\nabla_{e_i}v}u^{\flat}(e_i)\\
&=&-\nabla_{e_i}\left(g(\nabla_{e_i}u,v)\right)
+ \nabla_{e_i}\left(g(\nabla_vu,e_i)\right) +
g\left(\nabla_{\nabla_{e_i}e_i}u,v\right)
- g\left(\nabla_vu,\nabla_{e_i}e_i\right)\\
&&+ g\left(\nabla_{e_i}u,\nabla_{e_i}v\right)
- g\left(\nabla_{\nabla_{e_i}v}u,e_i\right)\\
&=&-g(\nabla_{e_i}\nabla_{e_i}u,v)-g(\nabla_{e_i}u,\nabla_{e_i}v)+g(\nabla_
{e_i}\nabla_vu,e_i)+g(\nabla_vu,\nabla_{e_i}e_i)\\
&&+ g\left(\nabla_{\nabla_{e_i}e_i}u,v\right)
- g\left(\nabla_vu,\nabla_{e_i}e_i\right )
+g\left(\nabla_{e_i}u,\nabla_{e_i}v\right)
-g\left(\nabla_{\nabla_{e_i}v}u,e_i\right)\\
&=& - g\left(\nabla^2_{e_i,e_i}u,v\right)
+g\left(\nabla_{e_i}\nabla_vu,e_i\right)
-g\left(\nabla_{\nabla_{e_i}v}u,e_i\right).
\end{eqnarray*}
Using the formula for the curvature $R $ and the Ricci
curvature we obtain
\begin{eqnarray*}
-(d\delta+\delta
d)u^\flat(v)&=&g(\nabla^2_{e_i,e_i}u,v)-g(\operatorname{R}
(e_i,v)u,e_i)+g(\nabla_{e_i}u,\nabla_ve_i)+g(\nabla_{\nabla_ve_i}u,e_i)\\
&=&g(\nabla^2_{e_i,e_i}u,v)-\operatorname{Ricci}(u,v)+0\\
&=&g(\nabla^2_{e_i,e_i}u-\operatorname{Ric}(u),v).
\end{eqnarray*}
The fact that $g(\nabla_{e_i}u,\nabla_ve_i)+g(\nabla_{\nabla_ve_i}u,e_i)=0$
can be simply proved pointwise at $x\in M$, assuming $\nabla e_i(x)=0$. (See
\cite{Pe1997} p.176/7 for a proof for a general local
orthonormal frame).$\,\,
\blacksquare$

\begin{lemma}
\label{divnabla}
For all $u,v\in\mathfrak{X}^{C^2}(M)$ we have
\[
\operatorname{div}(\nabla_vu)=\operatorname{Tr}(\nabla u\cdot\nabla v)
+\operatorname{Ricci}(u,v)+g(\operatorname{grad}\operatorname{div}(u),v)
\]
\end{lemma}
\textit{Proof} : We shall prove the identity at a fixed
point $x\in M$ so we  can choose a local orthonormal frame
$\{e_i\}$ such that $\nabla e_i(x)=0$. For the (1,1)
tensor $\nabla u $ we shall use the notation $\nabla u(v) :
= \nabla_v u$. At
$x$ we have :
\begin{align*}
\operatorname{Tr}(\nabla u\cdot\nabla v)
& +\operatorname{Ricci} (u,v)
=g\left(\nabla u \left(\nabla v(e_i)\right),e_i\right)
+ g \left(\operatorname{R}(e_i,v)u,e_i\right)\\
& = g\left(\nabla_{\nabla_{e_i}v}u,e_i\right)
+ g\left(\nabla_{e_i}\nabla_vu,e_i\right)
-g\left(\nabla_v\nabla_{e_i}u,e_i\right)
-g \left(\nabla_{[e_i,v]}u,e_i\right)\\
& =g(\nabla_{e_i}\nabla_vu,e_i)
-g(\nabla_v\nabla_{e_i}u,e_i) \quad \text{ because \quad
$\nabla_{e_i}v=[e_i,v]$ at $x$}\\
&=g(\nabla_{e_i}\nabla_vu,e_i)-\nabla_v(g(\nabla_{e_i}u,e_i))+g(\nabla_{e_i}
u,\nabla_ve_i)\\
&=\operatorname{div}(\nabla_vu)-d(\operatorname{div}(u))(v)+0\\
&=\operatorname{div}(\nabla_vu)-g(\operatorname{grad}\operatorname{div}(u),v)
\end{align*}
We can do that at each $x$ so the identity is proved.$\,\, \blacksquare$
\medskip

We shall denote below by $\Gamma^{L^2}(L(TM,TM))$ the
$L^2$ sections of the vector bundle $L(TM,TM)\mapsto M$.
\begin{lemma}
Consider on $\Gamma^{L^2}(L(TM,TM))$  the $L^2$ inner
product
\[
(R,S)_0:=\int_M\overline{g}(R,S)\mu.
\]
Then the following identities hold:\\
$(1)$ For all $u,v\in\mathfrak{X}^{C^2}(M)$ :
\[
(\nabla u,\nabla v)_0=\int_{\partial M}g(\nabla_nu,v)\mu_{\partial}-\langle
(\Delta+\operatorname{Ric})(u),v\rangle_0
\]\[
(\nabla u,(\nabla v)^t)_0=\int_{\partial M}g(\nabla_vu,n)\mu_{\partial}-
\langle (\operatorname{Ric}+\operatorname{grad}\operatorname{div})
(u),v\rangle_0.
\]
$(2)$ For all $u,v\in\mathfrak{X}_{|\!|}^{C^2}(M)$ :
\[
-2(\operatorname{Def}(u),\operatorname{Def}(v))_0=\langle\mathcal{L}
(u),v\rangle_0-\int_{\partial
M}g\left((\nabla_nu)^{tan}+S_n(u),v\right)\mu_{\partial}.
\]
Here $n$ denotes the outward-pointing unit normal vector
field along the boundary $\partial M$. We let $S_n :
T\partial M\rightarrow T\partial M$ be Weingarten map
defined by $S_n(u):=-\nabla_un$. The symbol
$(.)^{tan}$ denotes the tangential part to the  boundary of
a vector in $TM|\partial M$.
\end{lemma}
\textit{Proof} : (1) Let $\{e_i\}$ be a local orthonormal
frame on an open  subset $U$ of $M$. Recall the formula
$\operatorname{div}(f u) = f \operatorname{div}(u) + df(u)$.
On
$U$ we have :
\begin{eqnarray*}
\overline{g}(\nabla u,\nabla v)&=&\operatorname{Tr}((\nabla
u)^t\cdot\nabla v)
=g(\nabla_{e_i}u,\nabla_{e_i}v)\\
&=&d(g(\nabla_{e_i}u,v))(e_i)-g(\nabla_{e_i}\nabla_{e_i}u,v)\\
&=&\operatorname{div}(g(\nabla_{e_i}u,v)e_i)-g(\nabla_{e_i}u,v)\operatorname
{div}(e_i)-g(\nabla_{e_i}\nabla_{e_i}u,v).
\end{eqnarray*}
Using the relation $\nabla_{e_j} e_j = \sum_i g(\nabla_{e_j}
e_j, e_i) e_i $ in the third equality below, we get
\[
g(\nabla_{e_i}u,v)\operatorname{div}(e_i)=g(\nabla_{e_i}u,v)g(\nabla_
{e_j}e_i,e_j)=-g(\nabla_{e_i}u,v)g(e_i,\nabla_{e_j}e_j)
=g(\nabla_{\nabla_{e_j} e_j}u,v)
\]
and hence we conclude
\begin{align}
\label{inner product on nablas}
\overline{g}(\nabla u,\nabla
v)&=\operatorname{div}(g(\nabla_{e_i}u,v)e_i)-g
(\nabla^2_{e_i,e_i}u,v) \nonumber\\
&=\operatorname{div}(g(\nabla_{e_i}u,v)e_i)
-g((\Delta+\operatorname{Ric}) (u),v)
\end{align}
because of formula \eqref{Weitzenbock formula}.
We remark that the vector field $g(\nabla_{e_i}u,v)e_i$ does not depend on the
choice of the local orthonormal frame, so it defines a vector field on $M$.
Denote by $w$ this vector field. We obtain from
\eqref{inner product on nablas} using Stokes' theorem:
\begin{eqnarray*}
(\nabla u,\nabla v)_0&=&\int_M\overline{g}(\nabla u,\nabla v)
=\int_M\operatorname{div}(w)\mu-\int_Mg((\Delta+\operatorname{Ric})u,v)\mu\\
&=&\int_{\partial
M}g(w,n)\mu_{\partial}-\langle(\Delta+\operatorname{Ric})
(u),v\rangle_0.
\end{eqnarray*}
On $U$ we have $g(w,n)=g(g(\nabla_{e_i}u,v)e_i,n)
=g(\nabla_nu,v)$. So the
first identity is proved.

We proceed similarly with the proof of the second
identity. We have:
\begin{eqnarray*}
\overline{g}(\nabla u,(\nabla v)^t)&=&\operatorname{Tr}(\nabla
v\cdot\nabla u)
=g(e_i,\nabla_{\nabla_{e_i}u}v)\\
&=&d(g(e_i,v))(\nabla_{e_i}u)-g(\nabla_{\nabla_{e_i}u}e_i,v)\\
&=&\operatorname{div}(g(e_i,v)\nabla_{e_i}u)-g(e_i,v)\operatorname{div}(\nabla_
{e_i}u)-g(\nabla_{\nabla_{e_i}u}e_i,v).
\end{eqnarray*}
Using the formula $\operatorname{div}(\nabla_{e_i}u)
=\operatorname{Tr}(\nabla u\cdot\nabla e_i)
+\operatorname{Ricci}(u,e_i)
+g(\operatorname{grad}\operatorname{div}
(u),e_i)$ proved in Lemma \ref{divnabla},  we obtain
\begin{eqnarray*}
g(e_i,v)\operatorname{div}(\nabla_{e_i}u)&=&g(e_i,v)\operatorname{Tr}(\nabla
u\cdot\nabla e_i)+g(e_i,v)g((\operatorname{Ric}+\operatorname{grad}
\operatorname{div})(u),e_i))\\
&=&g(e_i,v)\operatorname{Tr}(\nabla e_i\cdot\nabla
u)+g((\operatorname{Ric}
+\operatorname{grad}\operatorname{div})(u),v)\\
&=&g(e_i,v)g(\nabla_{\nabla_{e_j}u}e_i,e_j)+g((\operatorname{Ric}+\operatorname
{grad}\operatorname{div})(u),v)\\
&=&-g(e_i,v)g(e_i,\nabla_{\nabla_{e_j}u}e_j)+g((\operatorname{Ric}
+\operatorname{grad}\operatorname{div})(u),v)\\
&=&-g(v,\nabla_{\nabla_{e_j}u}e_j)+g((\operatorname{Ric}+\operatorname{grad}
\operatorname{div})(u),v).
\end{eqnarray*}
Thus $\overline{g}(\nabla u,(\nabla v)^t)=\operatorname{div}(g(e_i,v)\nabla_
{e_i}u)-g((\operatorname{Ric}+\operatorname{grad}\operatorname{div})(u),v)$.
As before, the vector field $w:=g(e_i,v)\nabla_{e_i}u$ does not depend on the
choice of the local orthonormal frame. We obtain :
\begin{eqnarray*}
(\nabla u,(\nabla v)^t)_0&=&\int_M\overline{g}(\nabla u,(\nabla v)^t)
=\int_M\operatorname{div}(w)\mu-\int_Mg((\operatorname{Ric}+\operatorname{grad}
\operatorname{div})(u),v)\mu\\
&=&\int_{\partial M}g(w,n)\mu_{\partial}-\langle(\operatorname{Ric}
+\operatorname{grad}\operatorname{div})(u),v\rangle_0
\end{eqnarray*}
by Stokes' theorem. On $U$ we have
$g(w,n)=g(g(e_i,v)\nabla_{e_i}u,n)=g(\nabla_vu,n)$. So the
second identity is proved.\\
$(2)$ Using the two formulas in part (1) and the defintions
\[
\operatorname{Def}u=\frac{\nabla u+(\nabla u)^t}{2}
\quad \text{ and } \quad \mathcal{L}
=\Delta+2\operatorname{Ric}+\operatorname{grad}\operatorname{div}
\]
a direct computation gives
\[
-2(\operatorname{Def}u,\operatorname{Def}v)_0=\langle\mathcal{L}(u),v\rangle-
\int_{\partial M}g(\nabla_nu,v)\mu_{\partial}-\int_{\partial M}g(\nabla_vu,n)
\mu_{\partial}.
\]
If $u,v$ are tangent to the boundary, then on
$\partial M$ we get the relations
$g(\nabla_nu,v)=g((\nabla_nu)^{tan},v)$ and
$g(\nabla_vu,n)=d(g(u,n))(v)-g(u,\nabla_vn)=0+g(u,S_n(v))=g(S_n(u),v).\,\,
\blacksquare$

\medskip

Now we shall prove the following useful Lemma.

\begin{lemma}
\label{useful}
$(1)$ For $r\geq 1,\,\mathcal{L} :\mathfrak{X}^r\rightarrow\mathfrak{X}^{r-2}$
is a continuous linear map.\\
$(2)$ For all $u,v\in\mathcal{V}^r_{mix}$ with $r\geq
2$ we have
\[
\langle u,v\rangle_1=\langle
(1-\alpha^2\mathcal{L})u,v\rangle_0.
\]
\end{lemma}
\textit{Proof} : The first part is a direct verification.
To prove the second we use the preceding Lemma to obtain
$\langle u,v\rangle_1=\langle (1-\alpha^2
\mathcal{L})u,v\rangle_0$ for all $u,v\in\mathcal{V}_{mix}^{C^2}$. By the
Sobolev embedding theorem, the identity holds for all
$u,v\in\mathcal{V}_{mix}^s,  s>\frac{1}{2}dim M+2$. Using the
fact that $\mathcal{V}_{mix}^s$ is dense in
$\mathcal{V}_{mix}^2$ with the $H^2$ topology, and the fact
that
$\langle\,\,,\,\rangle_0$,$\langle\,\,,\,\rangle_1$, and
$\mathcal{L}$ are  continuous on $\mathfrak{X}^2$, the
identity holds for vector fields in
$\mathcal{V}^2_{mix}.\,\, \blacksquare$
\\

Using the previous lemma and solving a boundary value problem we can prove
(see \cite{Shkoller2000}) that for $r\geq 1$ the linear map
\[
(1-\alpha^2\mathcal{L}) : \mathcal{V}
_{mix}^r\longrightarrow\mathfrak{X}^{r-2}
\]
is a continuous isomorphism with inverse
\[
(1-\alpha^2\mathcal{L})^{-1} : \mathfrak{X}^{r-2}\longrightarrow
\mathcal{V}_{mix}^r.
\]
We recall from \cite{Shkoller2000} the two principal results
concerning the geometry of the  Lagrangian averaged Euler
equation (LAE-$\alpha$).

\begin{theorem}(Stokes decomposition)
\label{stokes decomposition}
For $r\geq 1 $ we have the following $\langle\,\,,\,\rangle_1-$orthogonal
decomposition:
\[\mathcal{V}_{mix}^r=\mathcal{V}^r_{\mu,mix}\oplus(1-\alpha^2\mathcal{L})^{-1}
\operatorname{grad}
H^{r-1}(M)\]
We denote by $\mathcal{P}_e :
\mathcal{V}^r_{mix}\longrightarrow\mathcal{V}^r_{\mu,mix}
$ the projection onto the first factor (\textit{Stokes projector}).\\
Then
\[\overline{\mathcal{P}} : T\mathcal{D}^s_{mix}|\mathcal{D}^s_{\mu,mix}
\longrightarrow T\mathcal{D}^s_{\mu,mix}\]
defined by $\overline{\mathcal{P}}(u_\eta):=[\mathcal
{P}_e(u_\eta\circ \eta^{-1})]\circ \eta$, is a $C^{\infty}$ bundle map.
\end{theorem}

\begin{theorem}
\label{geodesic theorem}
Let $\eta(t)\in \mathcal{D}^s_{\mu,D}$ be a curve in
$\mathcal{D}^s_{\mu,D}$ and let
$u(t):=TR_{\eta(t)^{-1}}(\dot\eta(t))=\dot\eta(t)\circ\eta(t)^{-1}
\in\mathcal{V}^s_{\mu,D}$.
Then the following properties are equivalent :\\
$(1)$ $\eta(t)$ is a geodesic of $(\mathcal{D}^s_{\mu,D},\mathcal{G}^1)$\\
$(2)$ $u(t)$ is a solution of LAE-$\alpha$ :
\[(1-\alpha^2\Delta_r)\partial_tu(t)+\nabla_{u(t)}[(1-\alpha^2\Delta_r)u(t)]-
\alpha^2\nabla u(t)^t\cdot\Delta_ru(t)=-\operatorname{grad}p(t)\]
$(3)$ $u(t)$ is a solution of :
\begin{equation}
\label{LAE}
\partial_tu(t)+\mathcal{P}_e\left(\nabla_{u(t)}u(t)
+\mathcal{F}^{\alpha}(u(t))\right)=0
\end{equation}
where $\mathcal{F}^{\alpha}:=\mathcal{U}^{\alpha}+\mathcal{R}^{\alpha} :
\mathcal{V}^s_{\mu,D}\longrightarrow\mathcal{V}^s_D$ with :
\begin{align}
\label{Ualpha}
\mathcal{U}^{\alpha}(u)&:=(1-\alpha^2\mathcal{L})^{-1}\alpha^2\operatorname
{Div}
(\nabla
u\cdot\nabla u^t+\nabla u\cdot\nabla u
-\nabla u^t\cdot\nabla u) \\
\label{Ralpha}
\mathcal{R}^{\alpha}(u)
&:=(1-\alpha^2\mathcal{L})^{-1}\alpha^2\Big{(}
\operatorname{Tr}\big{(}
\nabla_{\cdot}(\operatorname{R}(\cdot,u)u)
+\operatorname{R}(\cdot,u)\nabla_
{\cdot}u+\operatorname{R}(u,\nabla_{\cdot}u)\cdot\big{)}\nonumber\\
&\qquad -(\nabla_u\operatorname
{Ric})u-\nabla u^t\cdot\operatorname{Ric}(u)\Big{)}
\end{align}
$(4)$ $V(t):=\dot\eta(t)$ (Lagrangian velocity) is a
solution of :
\[\dot V(t)=\mathcal{S}^1(V(t))\]
where $\mathcal{S}^1\in\mathfrak{X}^{C^{\infty}}(T\mathcal{D}^s_{\mu,D})$ is
the geodesic spray of $(\mathcal{D}^s_{\mu,D},\mathcal{G}^1)$.
\end{theorem}

In part (3), $\operatorname{Div}$ denotes the divergence of a $(1,1)$-tensor :
\[
\operatorname{Div}(S):=(\nabla_{e_i}S)(e_i)
\]
for $\{e_i\}$ a local orthonormal frame.
In the last section we will generalise the previous theorem to the case of Neumann and mixed boundary conditions.
\section{Geodesic spray and connector of $(\mathcal{D}^s_{\mu,D},\mathcal{G}^1)
$}

In this section we shall give the formula of the geodesic spray $\mathcal{S}^1
$ and the connector $K^1$ of the weak Riemannian metric $\mathcal{G}^1$ on
$\mathcal{D}^s_{\mu,D}$. Recall that the geodesic spray is the Lagrangian
vector field on $T\mathcal{D}^s_{\mu,D}$ associated to the
Lagrangian $L :  T\mathcal{D}^s_{\mu,D}\longrightarrow
\mathbb{R}$ given by $L(u_{\eta})=\frac
{1}{2}\mathcal{G}^1(\eta)(u_\eta,u_\eta)$, that is
\[
\operatorname{\bf i}_{\mathcal{S}^1}\Omega_L=dL
\]
where $\Omega_L$ is the weak symplectic form associated to
$L$, that is, the pull back by the Legendre transformation
defined by $L $ of the canonical weak symplectic form on $T
^\ast \mathcal{D}^{s}_{\mu, D}$ (see, e.g.
\cite{MaRa1999}). So the integral curves of the geodesic
spray are
$V(t)=\dot\eta(t)
$ where $\eta(t)$ is a geodesic of $(\mathcal{D}^s_{\mu,D},\mathcal{G}^1)$.
Using that $u(t):=\dot\eta(t)\circ\eta(t)^{-1}$ is a solution of $\eqref{LAE}$
we will prove the following lemma.
\begin{lemma}
\label{spray}
The geodesic spray of $(\mathcal{D}^s_{\mu,D},\mathcal{G}^1)$ is given by :
\[\mathcal{S}^1(u_\eta)=T\overline{\mathcal{P}}\big{(}S\circ u_\eta
- \operatorname{Ver}_
{u_\eta}(\overline{\mathcal{F}}^{\alpha}(u_\eta))\big{)}\]
where $\overline{\mathcal{F}}^{\alpha}(u_\eta):=\mathcal{F}^{\alpha}
(u_\eta\circ \eta^{-1})\circ \eta$ and $S$ is the geodesic spray of
$(M,g)$ and $\operatorname{Ver}_{u_\eta}(v_\eta) \in
T_{u_\eta}(T\mathcal{D}^s_{\mu, D})$ is the vertical lift of
$v_\eta \in T_\eta
\mathcal{D}^s_{\mu, D}$ at $u_\eta \in T_\eta
\mathcal{D}^s_{\mu, D}$, that is,
\[
\operatorname{Ver}_{u_\eta}(v_\eta) =
\left.\frac{d}{dt}\right|_{t=0} ( u_\eta + t v_ \eta).
\]
\end{lemma}
\textit{Proof} : Let $\eta(t)$ be a goedesic of $(\mathcal{D}^s_
{\mu,D},\mathcal{G}^1)$. Then $u(t):=\dot\eta(t)\circ\eta(t)^{-1}$ is a
solution of
\[\partial_tu(t)+\mathcal{P}_e(\nabla_{u(t)}u(t)+\mathcal{F}^{\alpha}(u(t)))=0.
\]
We have $V(t)=\dot\eta(t)=u(t)\circ\eta(t)$. In the
following computation we denote by $\dot{u}(t)$ the
$t$-derivative  of $u(t) $ thought of as a curve in
$T\mathcal{D}^s_{\mu, D}$. However, $u(t) \in
\mathcal{V}^s_{\mu, D}$ for all $t $ and  therefore, one
can take the derivative $\partial_t u (t)$ of $u(t)$ as a
curve in the Hilbert space $\mathcal{V}^s_{\mu, D}$. The
relation between these two derivatives is $\dot{u}(t) =
\operatorname{Ver}_{u(t)}(\partial_t u (t))$ using the
standard identification between a vector space and its
tangent space at a point. Differentiating
$V(t)$ and using the  preceding equation we obtain
\begin{eqnarray*}
\dot V(t)&=&T(u(t))\circ\dot\eta(t)+\dot u(t)\circ\eta(t)\\
&=&T(u(t))\circ\dot\eta(t)+\operatorname{Ver}_{u(t)}(\partial_tu(t))\circ\eta
(t)\\
&=&T(u(t))\circ\dot\eta(t)-\operatorname{Ver}_{u(t)}(\mathcal{P}_e(\nabla_{u
(t)}u(t)+\mathcal
{F}^{\alpha}(u(t))))\circ\eta(t).
\end{eqnarray*}
We conclude that
\begin{eqnarray*}
\mathcal{S}^1(u_\eta)&=&T(u_\eta\circ \eta^{-1})\circ
u_\eta-\operatorname{Ver}_{u_\eta\circ
\eta^{-1}}\left(\mathcal{P}_e(\nabla_{u_\eta\circ
\eta^{-1}}(u_\eta\circ \eta^{-1})
+\mathcal{F}^{\alpha}(u_\eta\circ \eta^{-1}))\right)\circ
\eta\\
&=&Tu\circ u\circ \eta
-\operatorname{Ver}_u\big(\mathcal{P}_e(\nabla_u u
+\mathcal{F}^{\alpha}(u))\big)
\circ \eta \quad \text{ where $u:=u_\eta\circ \eta^{-1}
\in\mathcal{V}^s_{\mu,D}$}.
\end{eqnarray*}
Now it suffices to prove that for all
$u\in\mathcal{V}^s_{\mu,D}$ we have :
\begin{enumerate}

\item[(1)] $T\overline{\mathcal{P}}(Tu\circ u)
=Tu\circ u$ \quad and

\item[(2)]
$T\overline{\mathcal{P}}(Ver_u(\nabla_uu+\mathcal{F}^
{\alpha}(u)))
=Ver_u(\mathcal{P}_e(\nabla_uu+\mathcal{F}^{\alpha}(u)))$.

\end{enumerate}

\noindent(1) Let $c(t)$ be a curve in $\mathcal{D}^s_{\mu,D}$ such that $c(0)
=id_M$ and  $\dot c(0)=u$.
Let $d(t):=u\circ c(t)$. Then we have $d(0)=u$ and $\dot
d(0)=Tu\circ u$. We  get
\begin{eqnarray*}
T\overline{\mathcal{P}}(Tu\circ u)&=&
\left.\frac{d}{dt}\right|_{t=0}
\overline{\mathcal{P}}(u\circ c(t))
=\left.\frac{d}{dt}\right|_{t=0}\mathcal
{P}_e(u)\circ c(t)\\
&=&\left.\frac{d}{dt}\right|_{t=0}u\circ c(t)=Tu\circ u.
\end{eqnarray*}
(2) Let $v:=\nabla_uu+\mathcal{F}^{\alpha}(u)\in\mathcal{V}^{s-1}_D$. We get
\begin{eqnarray*}
T\overline{\mathcal{P}}(\operatorname{Ver}_u(v))
&=&T\overline{\mathcal{P}}\Big{(}\left.\frac{d}{dt}\right|_{t=0}
(u+tv)\Big{)}
=\left.\frac{d}{dt}\right|_{t=0}\overline
{\mathcal{P}}(u+tv)\\
&=&\left.\frac{d}{dt}\right|_{t=0}(u+t\mathcal{P}_e(v))=
\operatorname{Ver}_u(\mathcal{P}_e(v)).
\end{eqnarray*}
So, using  right-invariance  of $T\overline{\mathcal{P}}$ in
the second equality below and the expression of the spray
$S$ on $(M, g)$, namely $S \circ u = Tu \circ u -
\operatorname{Ver}_u(\nabla_u u)$, we obtain
\begin{eqnarray*}
\mathcal{S}^1(u_\eta)
&=&T\overline{\mathcal{P}}\big{(}Tu\circ u-
\operatorname{Ver}_u
(\nabla_uu+\mathcal{F}^{\alpha}(u))\big{)}\circ \eta\\
&=&T\overline{\mathcal{P}}\big{[}(Tu\circ u-
\operatorname{Ver}_u(\nabla_uu))\circ \eta
-\operatorname{Ver}_u
(\mathcal{F}^{\alpha}(u))\circ \eta\big{]}\\
&=&T\overline{\mathcal{P}}\big{(}S\circ u\circ\eta
-\operatorname{Ver}_{u_\eta}(\mathcal{F}^
{\alpha}(u)\circ \eta)\big{)}\\
&=&T\overline{\mathcal{P}}\big {(}S\circ
u_\eta-\operatorname{Ver}_{u_\eta}(\overline{\mathcal
{F}}^{\alpha}(u_\eta))\big{)}.\,\, \blacksquare
\end{eqnarray*}

Recall that locally the expressions of the geodesic spray
and the connector of a Riemannian manifold are given by
\[
\mathcal{S}^1(\eta,u)=(\eta,u,u,-\Gamma^1(\eta)(u,u))
\]and\[
\label{localconnector}
K^1(\eta,u,v,w)=(\eta,w+\Gamma^1(\eta)(u,v))
\]
where the symetric bilinear map $\Gamma^1(\eta)$ is the Christoffel map of the
Riemannian metric. Using these formula and the previous Lemma we
obtain  the global expression of $K^1$ below.

\begin{lemma}
\label{connector}
The connector $K^1 : TT\mathcal{D}^s_{\mu,D}\longrightarrow
T\mathcal{D}^s_{\mu,D}$ of
$(\mathcal{D}^s_{\mu,D},\mathcal{G}^1)$ is given by :
\[
K^1(X_{u_\eta})=\overline{\mathcal{P}}\Big{(}K\circ
X_{u_\eta}+\overline{\mathfrak
{F}}^{\alpha}\big{(}\pi_{_{T\mathcal{D}^s_{\mu,D}}}(X_{u_\eta}),T\pi_{_
{\mathcal{D}^s_{\mu,D}}}(X_{u_\eta})\big{)}\Big{)},
\]
where
\[
\displaystyle\overline{\mathfrak{F}}^{\alpha}(u_\eta,v_\eta):=\frac{1}
{2}\Big{(}\overline{\mathcal{F}}^{\alpha}(u_\eta+v_\eta)-\overline{\mathcal{F}}
^{\alpha}(u_\eta)-\overline{\mathcal{F}}^{\alpha}(v_\eta)\Big{)},
\]
$\pi_{_{\mathcal{D}^s_{\mu, D}}}: T\mathcal{D}^s_{\mu, D}
\rightarrow \mathcal{D}^s_{\mu, D}$ and
$\pi_{_{T\mathcal{D}^s_{\mu, D}}} : TT\mathcal{D}^s_{\mu, D}
\rightarrow T\mathcal{D}^s_{\mu, D}$ are tangent bundle
projections, and $K :  TTM\longrightarrow TM$ is the
connector of $(M,g)$.
\end{lemma}
\textit{Proof} : Let $\eta\in \mathcal{D}^s_{\mu,D}$,
$u_\eta,v_\eta\in  T_\eta\mathcal{D}^s_{\mu,D}$, and
$w_\eta\in T_\eta\mathcal{D}_D$. We  write
$\mathcal{S}^0(u_\eta):=S\circ u_\eta$ (in case $M$ has
no boundary, $\mathcal{S}^0$ is the geodesic spray of
$(\mathcal{D}^s,\mathcal{G}^0)$).

In local representation we have (with $(\eta,u), (\eta,v),
(\eta,w)$ the local  expressions of $u_\eta,v_\eta,w_\eta$):

\vspace{0,2cm}
$\mathcal{S}^i(\eta,u)=(\eta,u,u,-\Gamma^i(\eta)(u,u)), i=1,2$ where
$\Gamma^i$ are the Christoffel maps,

\vspace{0,2cm}
$\overline{\mathcal{F}}^{\alpha}(\eta,u)=(\eta,\overline{\mathcal{F}}^{\alpha}_
{loc}(\eta,u))$ and $\overline{\mathfrak{F}}^{\alpha}((\eta,u),(\eta,v))=
(\eta,\overline{\mathfrak{F}}^{\alpha}_{loc}(\eta)(u,v))$,

\vspace{0,2cm}
$\operatorname{Ver}_{(\eta,u)}(\overline{\mathcal{F}}^{\alpha}(\eta,u))=
(\eta,u,0,\overline
{\mathcal{F}}^{\alpha}_{loc}(\eta,u))$,

\vspace{0,2cm}
$\overline{\mathcal{P}}(\eta,u)=(\eta,\overline{\mathcal{P}}_{loc}(\eta,u))$,

\vspace{0,2cm}
$T\overline{\mathcal{P}}(\eta,u,v,w)=(\eta,\overline{\mathcal{P}}_{loc}
(\eta,u),v,D\overline{\mathcal{P}}_{loc}(\eta,u)(v,w))=(\eta,u,v,\overline
{\mathcal{P}}_{loc}(\eta,w))$.

\vspace{0,2cm}

\noindent Thus we find
\begin{eqnarray*}
\mathcal{S}^1(\eta,u)&=&T\overline{\mathcal{P}}(\mathcal{S}^0(\eta,u)
-\operatorname{Ver}_{(\eta,u)}(\overline{\mathcal{F}}^{\alpha}(\eta,u)))\text{
by Lemma \ref {spray}}\\
&=&T\overline{\mathcal{P}}(\eta,u,u,-\Gamma^0(\eta)(u,u)-\overline{\mathcal{F}}
^{\alpha}_{loc}(\eta,u))\\
&=&\left(\eta,u,u,-\overline{\mathcal{P}}_{loc}\left(\Gamma^0(\eta)(u,u)
+\overline
{\mathcal{F}}^{\alpha}_{loc}(\eta,u)\right)\right).
\end{eqnarray*}
We deduce that $\Gamma^1(\eta)(u,u)=\overline{\mathcal{P}}_{loc}
\left(\Gamma^0(\eta)
(u,u)+\overline{\mathcal{F}}^{\alpha}_{loc}(\eta,u)\right)$
and then that
\[
\Gamma^1(\eta)(u,v)=\overline{\mathcal{P}}_{loc}(\Gamma^0(\eta)(u,v)
+\overline{\mathfrak{F}}^{\alpha}_{loc}(\eta)(u,v)).
\]
Thus, with $u_\eta,v_\eta,w_\eta\in
T_\eta\mathcal{D}^s_{\mu,D}$, we obtain
\begin{eqnarray*}
K^1(\eta,u,v,w)&=&(\eta,w+\Gamma^1(\eta)(u,v))\\
&=&\left(\eta, w + \overline{\mathcal{P}}_{loc}
\left(\Gamma^0(\eta)(u,v)+\overline{\mathfrak
{F}}^{\alpha}_{loc}(\eta)(u,v)\right)\right)\\
&=&\overline{\mathcal{P}}\left((\eta,w+\Gamma^0(\eta)(u,v))
+(\eta,\overline
{\mathfrak{F}}^{\alpha}_{loc}(\eta)(u,v))\right)\\
&=&\overline{\mathcal{P}}\left(K^0(\eta,u,v,w)
+\overline{\mathfrak{F}}^{\alpha}
((\eta,u),(\eta,v))\right),
\end{eqnarray*}
where $K^0(X_{u_\eta}):=K\circ X_{u_\eta}$ with $K:TTM\longrightarrow TM$ the
connector of $(M,g)$. A globalisation of the previous formula gives the
result :
\[K^1(X_{u_\eta})=\overline{\mathcal{P}}\Big{(}K\circ X_{u_\eta}+\overline
{\mathfrak{F}}^{\alpha}\big{(}\pi_{_{T\mathcal{D}^s_{\mu,D}}}(X_{u_\eta}),T\pi_
{_{\mathcal{D}^s_{\mu,D}}}(X_{u_\eta})\big{)}\Big{)}.\,\, \blacksquare
\]

\section{The Lie-Poisson structure of LAE-$\alpha$ equation}
\label{section: l p structure}

In this section we shall define a Lie-Poisson bracket on
a certain class of functions on $\mathcal{V}_{\mu, D}^r$,
if $r \geq s > \frac{1}{2}\operatorname{dim}M +1 $ and
shall specify precise sharp conditions on their smoothness
class. In particular, we shall also determine the
conditions under which the Jacobi identity holds.

Let $s > \frac{1}{2}\operatorname{dim}M +1 $. Because of
the existence of the geodesic spray
$\mathcal{S}^1$ of the weak  Riemannian Hilbert manifold
$(\mathcal{D}^s_{\mu,D},\mathcal{G}^1)$ and the  fact that
the inclusion $\flat :
T_\eta\mathcal{D}^s_{\mu,D}\longrightarrow
T^*_\eta\mathcal{D}^s_{\mu,D}$ is dense, we can use the
results of section 4  in \cite{VaMa2004}. Therefore, by
those results, $T\mathcal{D}^s_{\mu,D}$  carries a Poisson
structure in the precise sense given there. To give it
explicitly in our case for the metric $\mathcal{G}^1$ we
need a few preliminaries.

If $F: T \mathcal{D}^s_{\mu,D} \rightarrow \mathbb{R}$ is
of class $C ^1$ we define the {\bfi horizontal
derivative\/} of $F $ by
\[
\frac{\partial F}{\partial \eta} : T
\mathcal{D}^s_{\mu,D}
\rightarrow T ^\ast \mathcal{D}^s_{\mu,D}
\]
by
\[
\left\langle \frac{\partial F}{\partial \eta}(u_\eta),
v_\eta  \right\rangle
: = \left.\frac{d}{dt}\right|_{t=0}F(\gamma(t)),
\]
where $\langle \,,\rangle$ is the duality paring
and $\gamma(t) \subset T \mathcal{D}^s_{\mu,D}$ is a smooth
path defined in a neighborhood of zero, with base point
denoted by $\eta(t) \subset
\mathcal{D}^s_{\mu,D}$, satisfying the following conditions:
\begin{itemize}
\item $\gamma(0) = u_\eta$
\item $\dot{\eta}(0) = v_\eta$
\item $\gamma$ is parallel, that is, its  covariant
derivative of the $\mathcal{G}^1$ Levi-Civita connection
vanishes.
\end{itemize}
The {\bfi vertical derivative\/}
\[
\frac{\partial F}{\partial u} : T
\mathcal{D}^s_{\mu,D}
\rightarrow T ^\ast \mathcal{D}^s_{\mu,D}
\]
of $F $ is defined as
the usual fiber derivative, that is,
\[
\left\langle \frac{\partial F}{\partial u}(u_\eta), v_\eta
\right\rangle : = \left.\frac{d}{dt}\right|_{t=0}
F(u_\eta + t v_\eta).
\]

These derivatives naturally induce corresponding functional
derivatives relative to the weak Riemannian metric
$\mathcal{G}^1$. The {\bfi horizontal\/} and {\bfi
vertical functional derivatives\/}
\[
\frac{\delta F}{\delta \eta }, \frac{\delta F}{\delta u} :
T \mathcal{D}^s_{\mu,D}
\rightarrow T \mathcal{D}^s_{\mu,D}
\]
are defined by the equalities
\[
\mathcal{G}^1(\eta)\left(\frac{\delta F}{\delta
\eta}(u_\eta),v_\eta \right)
= \left\langle\frac{\partial F}{\partial
\eta}(u_\eta), v_\eta \right\rangle
\quad \text{ and } \quad
\mathcal{G}^1 (\eta)\left(\frac{\delta F}{\delta
u}(u_\eta),v_\eta \right)
= \left\langle \frac{\partial F} {\partial u}(u_\eta),
v_\eta \right\rangle
\]
for any $u_\eta, v_ \eta\in T\mathcal{D}^s_{\mu,D}$. Note
that due to the weak character of $\mathcal{G}^1$, the
existence of the fucntional derivatives is not guaranteed.
But if they exist, they are unique.

We define, for $k\geq 1$ and
$r,t>\frac{1}{2}\operatorname{dim}M +1$ :
\[
C^k_r(T\mathcal{D}^t_{\mu,D}):=\left\{F\in
C^k(T\mathcal{D}^t_{\mu,D}) \Big|\exists\,\frac{\delta
F}{\delta \eta},\frac{\delta F}{\delta u} : T\mathcal{D}
^t_{\mu,D}\longrightarrow T\mathcal{D}^r_{\mu,D}\right\}.
\]

With these definitions the Poisson bracket of $F,G\in
C^k_r(T\mathcal{D}^t_{\mu,D})$ is given by
\begin{equation}
\label{poisson bracket one}
\{F,G\}^1(u_\eta)
=\mathcal{G}^1(\eta)\left(\frac{\delta
F}{\delta \eta} (u_\eta),\frac{\delta G}{\delta
u}(u_\eta)\right)
-\mathcal{G}^1(\eta)\left(
\frac{\delta F}{\delta u}(u_\eta), \frac{\delta G}{\delta
\eta}(u_\eta)\right)
\end{equation}
As in the case of Euler equation (see \cite{VaMa2004}) we
have the following  result.

\begin{proposition}
\label{commutative diagram}
Let $\pi_R : T\mathcal{D}^s_{\mu,D}\longrightarrow
\mathcal{V}^s_{\mu,D}$ be
definied by $\pi_R(u_\eta):=u_\eta\circ \eta^{-1}$.
Let $F_t$ be the flow of $\mathcal{S}^1$ and
$\widetilde{F}_t:=\pi_R\circ  F_t$.
Then $\widetilde{F}_t$ is the flow of LAE-$\alpha$ equation.
Moreover we have the  following commutative diagram :
$$\xymatrix{
T\mathcal{D}^s_{\mu,D} \ar[r]^{F_t} \ar[d]_{\pi_R}&
T\mathcal{D}^s_{\mu,D} \ar [d]^{\pi_R}\\
\mathcal{V}^s_{\mu,D}\ar[r]^{\widetilde{F}_t} &
\mathcal{V}^s_{\mu,D}.\\
}$$
\end{proposition}

\textit{Proof} : Let $u \in \mathcal{V}^s_{\mu,D}$ and
$V(t) = F_t (u)$. Then $V $ is an integral curve of
$\mathcal{S}^1$ with initial condition $u $. Note that
$\widetilde{F}_t (u) = \pi_R(V(t)) = V(t) \circ
\eta(t)^{-1}$, where $\eta(t) $ is the base point of
$V(t)$, which by Theorem \ref{geodesic theorem} (1) is the
geodesic of $\mathcal{S}^1$. Therefore, by Theorem
\ref{geodesic theorem} (2), $\widetilde{F}_t (u)$ is the
integral curve of LAE-$\alpha$ with initial condition $u$.

We still need to show that $\widetilde{F}_t \circ \pi_R =
\pi_R \circ F_t $. Indeed, since $ \mathcal{S}^1$ is a
right invariant vector field, its flow $F_t $ is right
equivariant and we conclude
\begin{align*}
(\widetilde{F}_t \circ \pi_R)(u_\eta)& = (\pi_R \circ F_t
\circ \pi_R)(u_\eta)
= (\pi_R \circ F_t \circ TR_{\eta^{-1}})(u_\eta) \\
&= (\pi_R \circ TR_{\eta^{-1}} \circ  F_t )(u_\eta)
= (\pi_R \circ  F_t )(u_\eta). \,\, \, \blacksquare
\end{align*}

We shall need later the fact that $\pi_R\in
C^k(T\mathcal{D}^{s+k}_{\mu,D},\mathcal{V}^s_{\mu,D})$ so
if $k = 0 $ then $\pi_R $ is only continuous.

\medskip

Our goal is to first  study the Lie-Poisson structure
of $\mathcal{V}^s_ {\mu,D}$ and secondly to show in what
sense the maps $F_t,\pi_R,\tilde{F}_t$  are Poisson maps.
We begin with the definition of some function spaces needed
later when we introduce the relevant Poisson bracket.

\begin{definition}  Let $s > \frac{1}{2}\operatorname{dim}M
+ 1 $.\\
$(1)$ For $k,t \geq 1$ and $ r\geq s$ define:
\[
C^k_{r,t}(\mathcal{V}^s_{\mu,D})
:=\{f\in C^k(\mathcal{V}^s_{\mu,D}) |\exists\,\delta f
:\mathcal{V}^r_{\mu,D}\longrightarrow
\mathcal{V}^t_{\mu,D}\} \quad
\text{ and }\quad
C^k_t(\mathcal{V}^s_{\mu,D})
:=C^k_{s,t}(\mathcal{V}^s_{\mu,D})
\]
where $\delta f$ is the {\bfi functional derivative\/} of f
with respect to  the inner product $\langle
\,\,,\,\rangle_1$:
\[
\langle\delta f(u),v\rangle_1=Df(u)(v),\quad
\forall\,u,v\in\mathcal{V}^r_{\mu,D}
\]

\noindent$(2)$ For $k\geq 0$, $r\geq s$, and $t\geq 1$
define:
\[\mathcal{K}^k_{r,t}(\mathcal{V}^s_{\mu,D}):=\{f\in C^{k+1}_{r,t}(\mathcal{V}
^s_{\mu,D})|\delta f\in C^k(\mathcal{V}^r_{\mu,D},\mathcal{V}^t_{\mu,D})\}\text
{ and }\mathcal{K}^k(\mathcal{V}^s_{\mu,D}):=\mathcal{K}^k_{s,s}(\mathcal{V}^s_
{\mu,D}).\]

\noindent$(3)$ Let $k\geq 1$, $r\geq s$, and
$t>\frac{1}{2}\operatorname{dim}M+1$. The
\textbf{Poisson
bracket on $C^k_{r,t}(\mathcal{V}^s_{\mu,D})$} is defined
by:
\begin{equation}
\label{Lie-Poisson bracket}
\{f,g\}^1_+(u):=\langle u,[\delta g(u),\delta
f(u)]\rangle_1,\quad \forall  u\in\mathcal{V}^r_{\mu,D}.
\end{equation}
\end{definition}

\noindent\textbf{Remark}
When $t>\frac{1}{2}\operatorname{dim} M+2$ we have
\[\{f,g\}^1_+(u)=\langle u,[\delta g(u),\delta f(u)]^R_{Lie}\rangle_1\]
where $[\,\,,\,]^R_{Lie}$ is the right-Lie bracket on the
``Lie-algebra" of  $\mathcal{D}^s_{\mu,D}$. We recognize
the classical Lie-Poisson bracket.
\medskip

Theorems \ref{properties1} and \ref{properties2} will
summarize the properties  of this Poisson bracket. In
the proofs we will use the three following Lemmas.

\begin{lemma} \label{1}
Let $s > \frac{1}{2}\operatorname{dim}M
+ 1 $.\\
$(1)$ Let $u\in \mathfrak{X}^s_{\operatorname{div},|\!|}$ and $v,w\in\mathfrak
{X}^s$. Then:
\[\langle v,\nabla_uw\rangle_0=-\langle\nabla_uv,w\rangle_0\]
$(2)$ Let $u,v\in\mathcal{V}^s_{\mu,D}$ and $w\in\mathcal{V}^s_D$. Then:
\[\langle v,\nabla_uw\rangle_1=-\langle\nabla_uv+\mathcal{D}^\alpha
(u,v),w\rangle_1\]
where $\mathcal{D}^\alpha : \mathcal{V}^s_{\mu,D}\times\mathcal{V}^s_{\mu,D}
\longrightarrow\mathcal{V}^s_D$ is the bilinear continuous map given by
\begin{align*}
\mathcal{D}^\alpha(u,v)&:=\alpha^2(1-\alpha^2\mathcal{L})^{-1}\Big{(}
\operatorname{Div}(\nabla v\cdot\nabla u^t+\nabla v\cdot\nabla u)\\
&\qquad
\qquad+\operatorname{Tr}\big{(}\nabla_{\cdot}(\operatorname{R}(\cdot,u)
v)+\operatorname{R}(\cdot,u)\nabla_\cdot v\big{)}\\
&\qquad \qquad+\operatorname{grad}\big{(}\operatorname{Tr}(\nabla
u\cdot\nabla
v)+\operatorname{Ricci}(u,v)\big{)}-(\nabla_u\operatorname{Ric})(v)\Big{)}
\end{align*}
\end{lemma}
\textit{Proof} : The first part follows by an integration
by parts argument which is justified since  all vector
fields are of class $C^1$ by the Sobolev embedding
theorem. Indeed, integrating the identity
$\pounds_u (g (v, w)) = \nabla_u (g (v, w)) = g (\nabla_u
v, w) + g (v, \nabla_u w)$ and using $\pounds_u \mu =
(\operatorname{div} u) \mu = 0 $ we get
\begin{align*}
&\langle \nabla_u v, w\rangle_0 + \langle v, \nabla_u
w\rangle_0
= \int_M g \left(\nabla_u v, w  \right)\mu +
\int_M g \left( v, \nabla_u w \right)\mu\\
&\quad = \int_M \pounds_u (g (v, w)) \mu =
 \int_M \pounds_u (g (v, w) \mu ) =
 \int_M d\mathbf{i}_u (g (v, w) \mu)\\
&\quad  = \int_{\partial M} \mathbf{i}_u (g (v, w) \mu) =
\int_{\partial M} g (v, w) g(u, n)\mu_\partial = 0
\end{align*}
by the Stokes theorem and the hypothesis that $g(u, n) = 0$
on $\partial M$.

For the second part we will use the following formula (see Lemma 3 in \cite
{Shkoller2000}): for all $u\in\mathcal{V}^s_{\mu,D}$ and $v\in\mathcal{V}^r_
{\mu,D},r>\frac{1}{2}\operatorname{dim}M+3$ we have:

\begin{equation}
\label{lemma3}
(1-\alpha^2\mathcal{L})^{-1}\nabla_u[(1-\alpha^2\Delta_r)v]=\nabla_uv+\mathcal
{D}^\alpha(u,v)
\end{equation}

Using Lemma \ref{useful}, the first part, and formula \eqref{lemma3} we obtain
for $u\in\mathcal{V}^s_{\mu,D}, w\in\mathcal{V}^s_D$ and $v\in\mathcal{V}^r_
{\mu,D},r>\frac{1}{2}\operatorname{dim}M+3$:
\begin{eqnarray*}
\langle v,\nabla_uw\rangle_1&=&\langle
(1-\alpha^2\Delta_r)v,\nabla_uw\rangle_0
\\
&=&-\langle\nabla_u[(1-\alpha^2\Delta_r)v],w\rangle_0\\
&=&-\langle (1-\alpha^2\mathcal{L})^{-1}\nabla_u[(1-\alpha^2\Delta_r)
v],w\rangle_1\\
&=&-\langle\nabla_uv+\mathcal{D}^{\alpha}(u,v),
w\rangle_1.
\end{eqnarray*}
Using the fact that
$v\in\mathcal{V}^r_{\mu,D},r>\frac{1}{2}\operatorname{dim}
M+3$ is dense in $\mathcal{V}^s_{\mu,D}$, and the fact that
$\langle\,,\rangle_1, \nabla$, and $\mathcal{D}^\alpha$ are continuous on
$\mathcal{V}^s_{\mu,D}$ we obtain that

\[\langle v,\nabla_uw\rangle_1=-\langle\nabla_uv+\mathcal{D}^{\alpha}
(u,v),w\rangle_1,\text{ for all }u,v\in\mathcal{V}^s_{\mu,D}\text{ and }
w\in\mathcal{V}^s_D.\,\, \blacksquare\]

\begin{lemma} \label{2}
Let $s>\frac{1}{2}\operatorname{dim} M+1$.
Let $B^\alpha : \mathcal{V}^{s+1}_{\mu,D}\times\mathfrak{X}
^s\longrightarrow\mathcal{V}^{s+1}_{\mu,D}$ the continuous bilinear map given
by
\[B^\alpha(v,w):=\mathcal{P}_e(1-\alpha^2\mathcal{L})^{-1}(\nabla w^t\cdot(1-
\alpha^2\Delta_r)v).\]
Then we have
\[\langle v,\nabla_uw\rangle_1=\langle B^\alpha(v,w),u\rangle_1\]
for all $u\in\mathcal{V}^r_{\mu,D},r>\frac{1}{2}\operatorname{dim}M$, and
for all $v\in\mathcal{V}^{s+1}_{\mu,D},$ and $w\in\mathfrak{X}^s$.
\end{lemma}
\textit{Proof} :
Using Lemma \ref{useful} and the Stokes decomposition
(see Theorem \ref{stokes decomposition}), we obtain:
\begin{eqnarray*}
\langle
v,\nabla_uw\rangle_1&=&\langle(1-\alpha^2\Delta_r)v,\nabla_uw\rangle_0
\\
&=&\langle\nabla w^t\cdot(1-\alpha^2\Delta_r)v,u\rangle_0\\
&=&\langle(1-\alpha^2\mathcal{L})^{-1}(\nabla
w^t\cdot(1-\alpha^2\Delta_r)
v),u\rangle_1\\
&=&\langle\mathcal{P}_e(1-\alpha^2\mathcal{L})^{-1}(\nabla
w^t\cdot(1-\alpha^2
\Delta_r)v),u\rangle_1. \;\; \blacksquare
\end{eqnarray*}

\begin{lemma}\label{3}
Let $s>\frac{1}{2}\operatorname{dim} M+1$.
Let $k\geq 1$ and $ f\in C^k(\mathcal{V}^s_{\mu,D})$
be such that there exists
$\delta f\in
C^1(\mathcal{V}^r_{\mu,D},\mathcal{V}^t_{\mu,D}),r\geq
s,t\geq 1$. Then:
\[\langle D\delta f(u)(v),w\rangle_1=\langle D\delta f(u)
(w),v\rangle_1,\forall\,u,v,w\in\mathcal{V}^r_{\mu,D}\]
\end{lemma}
\textit{Proof} : The proof is similar to that of Lemma
5.5 in \cite{VaMa2004}. We have
\begin{align*}
\langle D\delta f(u)(v),w\rangle_1
&= \left.\frac{d}{dt}\right|_{t=0} \langle \delta f (v +
tu ), w \rangle_1
 = \left.\frac{d}{dt}\right|_{t=0} Df (v + tu ) ( w ) \\
& = \left.\frac{d}{dt}\right|_{t=0}
\left.\frac{d}{ds}\right|_{s=0} f (v + tu + sw) \\
&= \left.\frac{d}{ds}\right|_{s=0}
\left.\frac{d}{dt}\right|_{t=0}f (v + tu + sw) \\
& = \langle D\delta f(u) (w),v\rangle_1.
\;\; \blacksquare
\end{align*}

\begin{theorem}
\label{properties1}
Let $s>\frac{1}{2}\operatorname{dim}M+1$ and $k\geq 1$. Then:
\[\{\,\,,\,\}^1_+:\mathcal{K}^k(\mathcal{V}^s_{\mu,D})\times\mathcal{K}^k
(\mathcal{V}^s_{\mu,D})\longrightarrow\mathcal{K}^{k-1}_{s+1,s-1}(\mathcal{V}
^s_{\mu,D})\]
and for all $u\in\mathcal{V}^{s+1}_{\mu,D}$ we have
\begin{align*}
&\delta(\{f,g\}^1_+)(u)=\mathcal{P}_e(\nabla_{\delta g(u)}\delta
f(u)-\nabla_
{\delta f(u)}\delta g(u))\\
&\qquad+D\delta g(u)\left(\mathcal{P}_e\left(\nabla_{\delta
f(u)}u+\mathcal{D}
^\alpha(\delta f(u),u)\right)+B^\alpha(u,\delta f(u)\right)\\
&\qquad-D\delta f(u)\left(\mathcal{P}_e\left(\nabla_{\delta
g(u)}u+\mathcal{D}
^\alpha(\delta g(u),u)\right)+B^\alpha(u,\delta g(u)\right)
\end{align*}
\end{theorem}
\textit{Proof} : Let $h:=\{f,g\}^1_+$. We have to show that $h\in\mathcal{K}^
{k-1}_{s+1,s-1}(\mathcal{V}^s_{\mu,D})$.\\
$\bullet$ Let's show that $h\in C^k(\mathcal{V}^s_{\mu,D})$.\\
We have $h(u)=\langle u,\nabla_{\delta g(u)}\delta f(u)\rangle_1-\langle
u,\nabla_{\delta f(u)}\delta g(u)\rangle_1$. Using the facts that $\nabla :
\mathcal{V}^s_{\mu,D}\times\mathcal{V}^s_{\mu,D}\longrightarrow\mathcal{V}^{s-
1}_D$ and $\langle\,\,,\,\rangle_1 : \mathcal{V}^{s-1}_D\times\mathcal{V}^{s-1}
_D\longrightarrow \mathbb{R}$ are bilinear continuous
maps, and that
$\delta f,\delta g\in C^k(\mathcal{V}^s_{\mu,D},\mathcal{V}^s_{\mu,D})$
by hypothesis, we  obtain the result.\\
$\bullet$ Let's show that $h \in
C^k(\mathcal{V}^s_{\mu,D})$ admits a functional
derivative  $\delta h\in C^{k-1}(\mathcal{V}^{s+1}_
{\mu,D},\mathcal{V}^{s-1}_{\mu,D})$.

 Let
$u,v\in\mathcal{V}^{s+1}_{\mu,D}$. Using Lemmas
\ref{2}, \ref{1}, and \ref{3}  we obtain:
\begin{align*}
Dh(u)(v)&=\langle v,\nabla_{\delta g(u)}\delta f(u)\rangle_1+\langle
u,\nabla_
{D\delta g(u)(v)}\delta f(u)\rangle_1\\
&\qquad+\langle u,\nabla_{\delta g(u)}D\delta f(u)(v)\rangle_1-
(f\leftrightarrow g)\\
&=\langle v,\mathcal{P}_e(\nabla_{\delta g(u)}\delta f(u))\rangle_1+\langle
B^\alpha(u,\delta f(u)),D\delta g(u)(v)\rangle_1\\
&\qquad-\langle\nabla_{\delta g(u)}u+\mathcal{D}^\alpha(\delta
g(u),u),D\delta
f(u)(v)\rangle_1-(f\leftrightarrow g)\\
&=\langle v,\mathcal{P}_e(\nabla_{\delta g(u)}\delta f(u))\rangle_1+\langle
D\delta g(u)(B^\alpha(u,\delta f(u))),v\rangle_1\\
&\qquad-\langle D\delta f(u)\left(\mathcal{P}_e\left(\nabla_{\delta g(u)}
u+\mathcal{D}^\alpha(\delta g(u),u)\right)\right),v\rangle_1-(f\leftrightarrow
g).
\end{align*}
Thus we conclude that the functional derivative exists
and equals
\begin{align*}
&\delta h(u)=\mathcal{P}_e(\nabla_{\delta g(u)}\delta
f(u)-\nabla_ {\delta f(u)}\delta g(u))\\
&\qquad+D\delta g(u)\left(\mathcal{P}_e\left(\nabla_{\delta
f(u)}u+\mathcal{D}
^\alpha(\delta f(u),u)\right)+B^\alpha(u,\delta f(u)\right)\\
&\qquad-D\delta f(u)\left(\mathcal{P}_e\left(\nabla_{\delta
g(u)}u+\mathcal{D}
^\alpha(\delta g(u),u)\right)+B^\alpha(u,\delta g(u)\right)
\end{align*}
A meticulous analysis show that $\delta h\in C^{k-1}(\mathcal{V}^{s+1}_
{\mu,D},\mathcal{V}^{s-1}_{\mu,D})$.$\,\, \blacksquare$
\medskip

With all these preparations we can now establish the
precise sense in which \eqref{Lie-Poisson bracket} is a
Lie-Poisson bracket.

\begin{theorem}
\label{properties2}
Let $s,t>\frac{1}{2}\operatorname{dim} M+1,\,r\geq s,$ and $k\geq 1$.\\
$(1)$ $\{\,\,,\,\}^1_+$ is $\mathbb{R}$-bilinear and anti-symmetric on $C^k_
{r,t}(\mathcal{V}^s_
{\mu,D})\times C^k_{r,t}(\mathcal{V}^s_{\mu,D})$.\\
$(2)$ $\{\,\,,\,\}^1_+$ is a derivation in each factor:
\[\{fg,h\}_+^1=\{f,h\}_+^1g+f\{g,h\}_+^1,\forall\,f,g,h\in C^k_{r,t}(\mathcal
{V}^s_{\mu,D}).\]
$(3)$ If $s>\frac{1}{2}\operatorname{dim} M+2$, $\{\,\,,\,\}^1_+$satisfies
the Jacobi identity:\\
For all $f,g,h\in\mathcal{K}^k(\mathcal{V}^s_{\mu,D})$ and $u\in \mathcal{V}^
{s+1}_{\mu,D}$ we have:
\[\{f,\{g,h\}^1_+\}^1_+(u)+\{g,\{h,f\}^1_+\}^1_+(u)+\{h,\{f,g\}^1_+\}^1_+(u)=0
\]
\end{theorem}
\textit{Proof}: $(1)$ This is obvious.\\
$(2)$ A direct computation, using Lemma \ref{1}, the
fact that for all  $f,g\in
C^k_{r,t}(\mathcal{V}^s_{\mu,D})$ we have $fg\in
C^k_{r,t}(\mathcal{V}^s_{\mu,D})$, and the relation
$\delta (fg)(u) =\delta f(u)g(u)+f(u)\delta g(u)$
proves the required identity.\\
$(3)$ Let $f,g,h\in\mathcal{K}^k(\mathcal{V}^s_{\mu,D})$,
and $u\in\mathcal {V}^{s+1}_{\mu,D}$.
By Theorem \ref{properties1} we obtain
$\{g,h\}^1_+\in\mathcal{K}^{k-1}_
{s+1,s-1}(\mathcal{V}^s_{\mu,D})\subset
C^k_{s+1,s-1}(\mathcal{V}^s_{\mu,D})$.  Since
$s-1>\frac{1}{2}\operatorname{dim} M+1$ we can compute
the expression
$\{f,\{g,h\}^1_+\}
^1_+(u)$. Using Lemmas \ref{2}, \ref{1}, and \ref{3} we
obtain:
\begin{align*}
&\{f,\{g,h\}^1_+\}^1_+(u)\\
&\qquad=\langle u,[\delta \{g,h\}^1_+(u),\delta f(u)]\rangle_1\\
&\qquad=\langle u,\nabla_{\delta \{g,h\}^1_+(u)}\delta
f(u)\rangle_1-\langle
u,\nabla_{\delta f(u)}\delta \{g,h\}^1_+(u)\rangle_1\\
&\qquad=\langle B^{\alpha}(u,\delta f(u)),\delta \{g,h\}^1_+(u)
\rangle_1+\langle\nabla_{\delta f(u)}u+\mathcal{D}^\alpha(\delta f
(u),u),\delta \{g,h\}^1_+(u)\rangle_1\\
&\qquad=\langle \delta \{g,h\}^1_+(u),
B^{\alpha}(u,\delta
f(u))+\mathcal{P}_e\left(\nabla_{\delta f(u)}u+
\mathcal{D}^\alpha(\delta f(u),u)\right)\rangle_1 \\
&\qquad=\langle \delta \{g,h\}^1_+(u),B_f\rangle_1,
\end{align*}
where we denoted, for convenience, $B _f: =
B^{\alpha}(u,\delta f(u))
+\mathcal{P}_e\left(\nabla_{\delta f(u)}u+
\mathcal{D}^\alpha(\delta f(u),u)\right) \in
\mathcal{V}_{\mu, D}^s$. Using the formula in Theorem
\ref{properties1} this equals
\begin{align*}
&\langle \mathcal{P}_e(\nabla_{\delta h(u)}\delta
g(u)-\nabla_{\delta g
(u)}\delta h(u)),B_f\rangle_1\\
&\qquad\qquad + \langle D\delta
h(u)\left(\mathcal{P}_e\left(\nabla_{\delta g
(u)}u+\mathcal{D}^\alpha(\delta g(u),u)\right)
+B^\alpha(u,\delta g(u))
\right),B_f\rangle_1\\
&\qquad\qquad - \langle D\delta g(u)
\left(\mathcal{P}_e\left(\nabla_{\delta h(u)}u
+\mathcal{D}^\alpha(\delta h(u),u)\right)
+B^\alpha(u,\delta h(u))
\right),B_f\rangle_1\\
&\qquad =\langle [\delta h(u),\delta g(u)],
B^{\alpha}(u,\delta f(u))+
\nabla_ {\delta f(u)}u+\mathcal{D}^\alpha(\delta
f(u),u)\rangle_1 +D_{hgf}-D_{ghf},
\end{align*}
where we denote
\[
D_{hgf} : = \langle D\delta
h(u)\left(\mathcal{P}_e\left(\nabla_{\delta g
(u)}u+\mathcal{D}^\alpha(\delta g(u),u)\right)
+B^\alpha(u,\delta g(u)
\right),B_f\rangle_1.
\]
Note that by Lemma \ref{3}, we have $D_{hgf} = D_{hfg}$.
Using Lemma \ref{2} and \ref{1} this equals
\begin{align*}
&\langle
\nabla_{[\delta h(u),\delta g(u)]}\delta f(u),u\rangle_1-
\langle\nabla_{\delta f(u)}[\delta h(u),\delta g(u)],u\rangle_1+D_{hgf}-D_{ghf}
\\
&\qquad=\langle[[\delta h(u),\delta g(u)],\delta
f(u)],u\rangle_1+D_{hgf}-D_{ghf}
\\
&\qquad=\langle[[\delta h(u),\delta g(u)],\delta
f(u)],u\rangle_1+D_{hgf}-D_{gfh}.
\end{align*}
Using Jacobi identity for the Jacobi-Lie bracket of
vector fields we obtain:
\begin{align*}
&\{f,\{g,h\}^1_+\}^1_+(u)+\{g,\{h,f\}^1_+\}^1_+(u)+\{h,\{f,g\}^1_+\}^1_+(u)\\
&\qquad=0+(D_{hgf}-D_{gfh})+(D_{fhg}-D_{hgf})
+(D_{gfh}-D_{fhg}) =0 \;\; \blacksquare
\end{align*}

\section{Geometric Properties of the Flow of LAE-$\alpha$}

Now we will prove that the maps $\pi_R, F_t$, and
$\tilde{F}_t$ in Proposition \ref{commutative diagram}
are Poisson maps. As we shall see, the considerations
below need the hypothesis that  $\pi_R$ be at least of
class $C^1$. Note  that
$\pi_R:T\mathcal{D}^s_{\mu,D}\longrightarrow
\mathcal{V}^s_ {\mu,D}$ is only continuous. Later on we
shall use the fact that
$\pi_R\in C^k(T\mathcal{D}^{s+k}_{\mu,D},\mathcal{V}^s_
{\mu,D})$ for all $k \geq 0 $ (see \cite{EbMa1970}). If
$f\in C^k(\mathcal{V}^s_{\mu,D})$, we shall denote
$f_R:=f\circ\pi_R\in C^k (T\mathcal{D}^{s+k}_{\mu,D})$.

\begin{lemma}\label{vertical}
Let $k\geq 1$ and $r>\frac{1}{2}\operatorname{dim}M+1$
such that $s+k\geq r$. Let $f\in C^k_r
(\mathcal{V}^s_{\mu,D})$.
Then the vertical functional derivative of $f_R$ with respect to $\mathcal{G}^1
$ exists and is given by:
\[\frac{\delta f_R}{\delta u}(u_\eta)=TR_\eta(\delta f(\pi_R(u_\eta)))\in
T\mathcal{D}^r_{\mu,D}, \qquad \forall\,u_\eta\in
T\mathcal{D}^{s+k}_{\mu,D}\]
\end{lemma}
\textit{Proof}: This is a direct computation using the chain
rule, the right-invariance of $\mathcal{G}^1$, and the fact that the naturel
isomorphism
between a vector space and its tangent space at a point is the vertical-lift.
Indeed, we have:
\begin{align*}
\left\langle \frac{\partial f_R}{\partial u}
(u_\eta),v_\eta\right\rangle
&=\left.\frac{d}{dt}\right|_{t=0}f_R
(u_\eta+tv_\eta)
=\left.\frac{d}{dt}\right|_{t=0}
(f\circ\pi_R)(u_\eta+tv_\eta)\\
& = df(\pi_R(u_\eta)) \left(\left.\frac{d}{dt}\right|_{t=0}
\Big(\pi_R(u_\eta) + t \pi_R(v _\eta)
\Big)\right) \\
&=df(\pi_R(u_\eta))\left(\operatorname{Ver}_{\pi_R(u_\eta)}(\pi_R(v_\eta))
\right)
=Df(\pi_R(u_\eta)) (\pi_R(v_\eta))\\
&=\langle\delta
f(\pi_R(v_\eta)),\pi_R(v_\eta)\rangle_1
=\mathcal{G}^1(\eta)\Big
{(}TR_\eta(\delta f(\pi_R(v_\eta)),v_\eta\Big{)},
\end{align*}
where in the fifth equlity $D$ denotes the Fr\'echet
derivative of $f$ thought of as a function defined on the
Hilbert space $\mathcal{V}_{\mu,D}^s$ and in the third
equality $d$ denotes  the exterior derivative of $f $
thought of as a function defined on the manifold
$\mathcal{V}_{\mu,D}^s$.

So we conclude that the functional vertical covariant derivative exists and is
given by
\[
\frac{\delta f_R}{\delta u}(u_\eta)=TR_\eta(\delta
f(\pi_R(u_\eta))).
\]
Since $s+k\geq r$, it is an element of
$T\mathcal{D}^r_ {\mu,D}.\,\,\blacksquare$

\medskip

The computation of the horizontal functional derivative
of $f_R$ will  involve the connector and therefore the map
$\mathcal{F}^\alpha$ defined in Theorem \ref{geodesic
theorem}. The following Lemma  gives a useful expression
for $\mathcal{F}^\alpha$.

\begin{lemma}
\label{Falpha}
$(1)$ For all $u\in\mathcal{V}^s_{\mu,D}$ we have:
\begin{align*}\nabla u^t\cdot\Delta_r u&=\operatorname{Div}(\nabla
u^t\cdot\nabla u)-\operatorname{Tr}(\operatorname{R}(u,\nabla_\cdot u)\cdot)\\
&\qquad+\nabla u^t\cdot\operatorname{Ric}u-\frac{1}{2}\operatorname{grad}
(\operatorname{Tr}(g(\nabla_\cdot u,\nabla_\cdot u))).
\end{align*}
This shows that $\nabla u^t\cdot\Delta_r u$ is in
$\mathfrak{X}^{s-2}$.\\
$(2)$ For all $u\in\mathcal{V}^s_{\mu,D}$ we have:
\[
\mathcal{F}^{\alpha}(u)=\mathcal{D}^\alpha
(u,u)-(1-\alpha^2\mathcal{L})^{-1}
\alpha^2\Big{(}\operatorname{grad}(F(u))
+\nabla u^t\cdot\Delta_r u\Big{)},
\]
where $\mathcal{D}^\alpha$ was defined in Lemma \ref{1} and
$F\in C^\infty(\mathcal{V}^s_{\mu,D})$ is
given by
\[
F(u)=\operatorname{Tr}(\nabla u\cdot\nabla
u)+\operatorname{Ricci}(u,u)+\frac
{1}{2}\operatorname{Tr}(g(\nabla_\cdot u,\nabla_\cdot u)).
\]
\end{lemma}
\label{b}
\textit{Proof}:
(1) We shall prove the identity at a given point $x\in M$ so
we can choose a  local orthonormal frame $\{e_i\}$ such that
$\nabla e_i(x)=0$. The computation below is carried out at
the point $x $ and we shall not write this evaluation. We
have
\begin{align*}
&\operatorname{Div}(\nabla u^t\cdot\nabla u)=\nabla_{e_i}(\nabla
u^t\cdot\nabla u)(e_i)=\nabla_{e_i}(\nabla u^t\cdot\nabla u(e_i))\\
&\qquad=\nabla_{e_i}(g(\nabla u^t\cdot\nabla u(e_i),e_k)e_k)=\nabla_{e_i}(g
(\nabla_{e_i}u,\nabla_{e_k}u))e_k\\
&\qquad=g(\nabla_{e_i}\nabla_{e_i}u,\nabla_{e_k}u)e_k+g(\nabla_{e_i}u,\nabla_
{e_i}\nabla_{e_k}u)e_k\\
&\qquad=g(\nabla u^t\cdot\nabla_{e_i}\nabla_{e_i}u,e_k)e_k+g(\nabla_{e_i}
u,\operatorname{R}(e_i,e_k)u)e_k+g(\nabla_{e_i}u,\nabla_{e_k}\nabla_{e_i}u)
e_k\\
&\qquad=\nabla
u^t\cdot\nabla_{e_i}\nabla_{e_i}u+g(\operatorname{R}(u,\nabla_
{e_i}u)e_i,e_k)e_k+\frac{1}{2}\nabla_{e_k}(g(\nabla_{e_i}u,\nabla_{e_i}u))e_k\\
&\qquad=\nabla u^t\cdot\nabla_{e_i}\nabla_{e_i}u+\operatorname{R}(u,\nabla_
{e_i}u)e_i+\frac{1}{2}d(g(\nabla_{e_i}u,\nabla_{e_i}u))(e_k)e_k\\
&\qquad=\nabla u^t\cdot\nabla_{e_i}\nabla_{e_i}u+\operatorname{Tr}
(\operatorname{R}(u,\nabla_{\cdot}u)\cdot)+\frac{1}{2}g(\operatorname{grad}(g
(\nabla_{e_i}u,\nabla_{e_i}u)),e_k)e_k\\
&\qquad=\nabla u^t\cdot\nabla_{e_i}\nabla_{e_i}u+\operatorname{Tr}
(\operatorname{R}(u,\nabla_{\cdot}u)\cdot)+\frac{1}{2}\operatorname{grad}(g
(\nabla_{e_i}u,\nabla_{e_i}u))
\end{align*}
which, using  the Weitzenb\"ock formula in Lemma
\ref{Weitzenbock}, proves the desired formula.

(2) Using the
formulas (\ref{Ualpha}) and (\ref{Ralpha}), and part (1)
above, we have:
\begin{align*}
\mathcal{F}^\alpha(u)&=\mathcal{U}^\alpha(u)+\mathcal{R}^\alpha(u)\\
&=\mathcal{D}^\alpha(u,u)+(1-\alpha^2\mathcal{L})^{-1}\alpha^2\Big{(}-
\operatorname{grad}\left(\operatorname{Tr}(\nabla u\cdot\nabla u)+\operatorname
{Ricci}(u,u)\right)\\
&\qquad\qquad-\operatorname{Div}(\nabla u^t\cdot\nabla u)+\operatorname{Tr}
(\operatorname{R}(u,\nabla_{\cdot}u)\cdot)-\nabla u^t\cdot\operatorname{Ric}
u\Big{)}\\
&=\mathcal{D}^\alpha(u,u)-(1-\alpha^2\mathcal{L})^{-1}\alpha^2\Big{(}\nabla
u^t\cdot\Delta_r u\\
&\qquad\qquad +\operatorname{grad}\big{[}\operatorname{Tr}(\nabla
u\cdot\nabla
u)+\operatorname{Ricci}(u,u)+\frac{1}{2}\operatorname{Tr}(g(\nabla_\cdot
u,\nabla_\cdot u))\big{]}\Big{)}.\,\, \blacksquare
\end{align*}

\begin{lemma}
\label{horizontal}
Let $k\geq 1$ and $r>\frac{1}{2}\operatorname{dim}M+2$ such that $s+k\geq
r$.
Let $f\in C^k_r (\mathcal{V}^s_{\mu,D})$.
Then the horizontal functional derivative of $f_R$ with respect to $\mathcal{G}
^1$ exists. It is given by:
\[\frac{\delta f_R}{\delta \eta}(u_\eta)=\frac{1}{2}TR_h\Big{[}B^\alpha
(u,\delta f(u))-B^\alpha(\delta f(u),u)+\mathcal{P}_e\big{(}\mathcal{D}^\alpha
(\delta f(u),u)-\mathcal{D}^\alpha(u,\delta f(u))\big{)}\Big{]}\]
for all $u_\eta\in T\mathcal{D}^{s+k}_{\mu,D}$, where $u:=\pi_R(u_\eta)$
and $B ^\alpha$ was defined in Lemma \ref{2}. So
we have:
\[
\frac{\delta f_R}{\delta \eta}(u_\eta)\in
T\mathcal{D}^r_{\mu,D}, \qquad \forall\,u_\eta\in
T\mathcal{D}^{s+k}_{\mu,D}.
\]

\end{lemma}
\textit{Proof} : By Lemma \ref{connector}, we will have the
two following formula
\begin{equation}
\label{k1}
K^1(Tu\circ v)=\mathcal{P}_e(\nabla_vu
+\mathfrak{F}^\alpha(u,v)),
\end{equation}
where, using part (2) in Lemma \ref{Falpha} and the
definition of $\mathfrak{F}^\alpha$ in Lemma
\ref{connector}, we have
\begin{align}
\label{F-alpha}
\mathfrak{F}^\alpha(u,v)&=\frac{1}{2}\Big{(}\mathcal{F}^{\alpha}(u+v)-\mathcal
{F}^{\alpha}(u)-\mathcal
{F}^{\alpha}(v)\Big{)}\nonumber\\
&=\frac{1}{2}\Big{(}\mathcal{D}^\alpha(u,v)+\mathcal{D}^\alpha(v,u)\nonumber\\
&\qquad\qquad-(1-\alpha^2\mathcal{L})^{-1}\alpha^2\big{(}\operatorname{grad}(G
(u,v))+\nabla u^t\cdot\Delta_r v+\nabla v^t\cdot\Delta_r u\big{)}\Big{)}
\end{align}
denoting $G(u,v):=F(u+v)-F(u)-F(v)$.

Let $u_\eta,v_\eta\in T\mathcal{D}^{s+k}_{\mu,D}$ and $\gamma(t) \subset T
\mathcal{D}^s_{\mu,D}$ a smooth
path defined in a neighborhood of zero, with base point
denoted by $\eta(t) \subset
\mathcal{D}^s_{\mu,D}$, satisfying the following conditions:
\begin{itemize}
\item $\gamma(0) = u_\eta$
\item $\dot{\eta}(0) = v_\eta$
\item $\gamma$ is parallel.
\end{itemize}
By definition we have:
\begin{align*}
\left\langle\frac{\partial f_R}{\partial
\eta}(u_\eta),v_\eta\right\rangle&=\left.\frac{d}{dt}\right|_{t=0}(f\circ\pi_R)
(\gamma(t))\\
&=df(\pi_R(u_\eta))\left(\left.\frac{d}{dt}\right|_{t=0}\pi_R(\gamma(t))\right)
\\
&=df(\pi_R(u_\eta))\left(\operatorname{Ver}_{\pi_R(u_\eta)}\left(K^1\left
(\left.\frac{d}{dt}
\right|_{t=0}\pi_R(\gamma(t))\right)\right)\right)\\
&=Df(\pi_R(u_\eta))\left(K^1\left(\left.\frac{d}{dt}\right|_{t=0}\pi_R(\gamma
(t))\right)\right).
\end{align*}
For the third equality, it suffices to remark that $\left.\frac{d}{dt}\right|_
{t=0}\pi_R(\gamma(t))$ is a vertical vector field. To obtain the last equality
it suffices to use the natural isomorphism between a vector space and its
tangent space at a point.

Using the formulas for the derivative of the composition
and inversion we have:
\begin{align*}
\left.\frac{d}{dt}\right|_{t=0}\pi_R(\gamma(t))
&=\left.\frac{d}{dt}\right|_
{t=0}\left(\gamma(t)\circ\eta(t)^{-1}\right)\\
&=\left(\left.\frac{d}{dt}\right|_{t=0}\gamma(t)\right)
\circ\eta(0)^{-1}
+ T(\gamma(0))\circ
\left(\left.\frac{d}{dt}\right|_{t=0}\eta(t)^{-1}\right)\\
&=\left(\left.\frac{d}{dt}\right|_{t=0}\gamma(t)\right)
\circ\eta^{-1}
- Tu_\eta\circ T(\eta(0)^{-1})\circ
\left(\left.\frac{d}{dt}\right|_{t=0}\eta(t)\right)
\circ\eta(0)^{-1}\\
&=\left(\left.\frac{d}{dt}\right|_{t=0}\gamma(t)\right)
\circ\eta^{-1}-Tu_\eta\circ
T\eta^{- 1}\circ v_\eta\circ\eta^{-1}\\
&=\left(\left.\frac{d}{dt}\right|_{t=0}\gamma(t)\right)
\circ\eta^{-1}
-T(u_\eta\circ \eta^{-1}) \circ(v_\eta\circ\eta^{-1}).
\end{align*}
Using the right-invariance of the connector, the definition of the
covariant derivative $ D/dt$, and the fact that $\gamma(t)$
is parallel we obtain:
\begin{align*}
K^1\left(\left.\frac{d}{dt}\right|_{t=0}\pi_R(\gamma(t))\right)
&=K^1\left(\left(\left.\frac{d}{dt}\right|_{t=0}\gamma(t)\right)
\circ\eta^{-1}\right)
-K^1(T
(u_\eta\circ \eta^{-1})\circ (v_\eta\circ\eta^{-1}))\\
&=K^1\left(\left.\frac{d}{dt}\right|_{t=0}\gamma(t)\right)
\circ\eta^{-1}-K^1(T
(u_\eta\circ \eta^{-1})\circ (v_\eta\circ\eta^{-1}))\\
&=\left(\left.\frac{D}{dt}\right|_{t=0}\gamma(t)\right)
\circ\eta^{-1}-K^1(T(u_\eta\circ
\eta^{-1})\circ (v_\eta\circ\eta^{-1}))\\
&=0-K^1(T(u_\eta\circ \eta^{-1})\circ
(v_\eta\circ\eta^{-1})).
\end{align*}
Thus we obtain:
\begin{align*}
&\left\langle\frac{\partial f_R}{\partial
\eta}(u_\eta), v_\eta\right\rangle
=-Df(u_\eta\circ\eta^{-1})
(K^1(T(u_\eta\circ \eta^{-1})\circ(v_\eta\circ \eta^{-1})))\\
&\qquad =-Df(u)(K^1(Tu\circ v))\qquad \text{where
$u:=u_\eta\circ
\eta^{-1}$ and
$v:=v_\eta\circ \eta^{-1}$}\\
&\qquad
=-Df(u)(\mathcal{P}_e(\nabla_vu+\mathfrak{F}^\alpha(u,v)))
\qquad \text{by formula \eqref{k1}}\\
&\qquad =-\langle\delta
f(u),\mathcal{P}_e(\nabla_vu+\mathfrak{F}^\alpha(u,v))
\rangle_1\\
&\qquad =-\langle\delta
f(u),\nabla_vu+\mathfrak{F}^\alpha(u,v)\rangle_1\\
&\qquad =-\langle\delta
f(u),\nabla_vu+\frac{1}{2}(\mathcal{D}^\alpha(u,v)
+\mathcal{D}^\alpha(v,u))\rangle_1\\
&\qquad \qquad
+\frac{1}{2}\langle\delta f(u),(1-\alpha^2\mathcal{L})^{-1}
\alpha^2\big{(}\operatorname{grad}(G(u,v))\big{)}\rangle_1\\
&\qquad \qquad +\frac{1}{2}\langle\delta
f(u),(1-\alpha^2\mathcal{L})^{-1}
\alpha^2\big{(}\nabla u^t\cdot\Delta_r v+\nabla v^t\cdot\Delta_r u\big{)}
\rangle_1 \qquad \text{by formula \eqref{F-alpha}}.
\end{align*}
The second term is zero because of the Stokes decomposition
(see Theorem \ref{stokes decomposition}). For the first term
we  have by Lemmas \ref{1} and \ref{2}:
\begin{align*}
-\langle\delta f(u),&\nabla_vu
+\frac{1}{2}(\mathcal{D}^\alpha(u,v)
+\mathcal{D}^\alpha(v,u))\rangle_1\\
&=-\langle\delta
f(u),\frac{1}{2}\nabla_vu+\frac{1}{2}\mathcal{D}^\alpha(u,v)
+\frac{1}{2}(\mathcal{D}^\alpha(v,u)+\nabla_vu)\rangle_1\\
&=-\frac{1}{2}\Big{(}\langle\delta f(u),\nabla_vu\rangle_1+\langle\delta f
(u),\mathcal{D}^\alpha(u,v)\rangle_1-\langle\nabla_v\delta f(u),u\rangle_1\Big
{)} \\
&=-\frac{1}{2}\Big{(}\langle B^\alpha(\delta
f(u),u),v\rangle_1+\langle\delta f
(u),\mathcal{D}^\alpha(u,v)+\nabla_uv\rangle_1\\
&\qquad\qquad-\langle\delta f(u),\nabla_uv\rangle_1-\langle
B^\alpha(u,\delta f
(u)),v\rangle_1\Big{)}\\
&=-\frac{1}{2}\Big{(}\langle B^\alpha(\delta f(u),u),v\rangle_1-
\langle\nabla_u\delta f(u),v\rangle_1\\
&\qquad\qquad+\langle\nabla_u\delta f(u)+\mathcal{D}^\alpha(u,\delta f
(u)),v\rangle_1-\langle B^\alpha(u,\delta f(u)),v\rangle_1\Big{)}\\
&=\frac{1}{2}\langle B^\alpha(u,\delta f(u))-B^\alpha(\delta
f(u),u)-\mathcal
{D}^\alpha(u,\delta f(u)),v\rangle_1.
\end{align*}
By Lemmas \ref{useful}, \ref{1}, the third term becomes:
\begin{align*}
\langle\delta f(u),&(1-\alpha^2\mathcal{L})^{-1}\alpha^2\big{(}\nabla
u^t\cdot\Delta_r v+\nabla v^t\cdot\Delta_r u\big{)}\rangle_1\\
&=\langle\delta f(u),\alpha^2\big{(}\nabla u^t\cdot\Delta_r v+\nabla
v^t\cdot\Delta_r u\big{)}\rangle_0\\
&=-\langle\delta f(u),\nabla u^t\cdot (1-\alpha^2\Delta_r)
v\rangle_0+\langle\delta f(u),\nabla u^t\cdot v\rangle_0\\
&\qquad-\langle\delta f(u),\nabla v^t\cdot (1-\alpha^2\Delta_r)
u\rangle_0+\langle\delta f(u),\nabla v^t\cdot u\rangle_0\\
&=-\langle \nabla_{\delta f(u)}u,v\rangle_1+\langle \nabla_{\delta f(u)}
u,v\rangle_0\\
&\qquad-\langle \nabla_{\delta f(u)}v,u\rangle_1+\langle \nabla_{\delta
f(u)}
v,u\rangle_0\\
&=-\langle \nabla_{\delta f(u)}u,v\rangle_1+\langle \nabla_{\delta f(u)}
u,v\rangle_0\\
&\qquad+\langle\nabla_{\delta f(u)}u+\mathcal{D}^\alpha(\delta f
(u),u),v\rangle_1-\langle v,\nabla_{\delta f(u)}u\rangle_0\\
&=\langle\mathcal{D}^\alpha(\delta f(u),u),v\rangle_1.
\end{align*}
So we obtain:
\begin{align*}
&\frac{\partial f_R}{\partial \eta}(u_\eta)(v_\eta)\\
&=\frac{1}{2}\langle B^\alpha(u,\delta
f(u))-B^\alpha(\delta f(u),u)
+\mathcal{P}_e\big{(}\mathcal{D}^\alpha(\delta
f(u),u)-\mathcal{D}^\alpha (u,\delta
f(u))\big{)},v\rangle_1\\
&=\frac{1}{2}\mathcal{G}^1(\eta)\Big{(}TR_\eta\Big{[}B^\alpha(u,\delta
f (u))-B^\alpha(\delta
f(u),u)+\mathcal{P}_e\big{(}\mathcal{D}^\alpha(\delta f
(u),u)-\mathcal{D}^\alpha(u,\delta
f(u))\big{)}\Big{]},v_\eta\Big{)}.
\end{align*}
Therefore we obtain the existence of
$\displaystyle\frac{\delta f_R}{\delta
\eta}(u_\eta)\in T\mathcal{D}^r_{\mu,D}$, given by
\[\frac{\delta f_R}{\delta \eta}(u_\eta)=\frac{1}{2}TR_h\Big{[}B^\alpha
(u,\delta f(u))-B^\alpha(\delta f(u),u)+\mathcal{P}_e\big{(}\mathcal{D}^\alpha
(\delta f(u),u)-\mathcal{D}^\alpha(u,\delta
f(u))\big{)}\Big{]},\] where
$u:=\pi_R(u_\eta).\,\,\blacksquare$\\

Lemmas \ref{vertical} and \ref{horizontal} yield the
following theorem:

\begin{theorem}\label{f_R}
Let $k\geq 1$ and $r>\frac{1}{2}\operatorname{dim}M+2$ such
that $s+k\geq r$. Let $f\in C^k_r
(\mathcal{V}^s_{\mu,D})$. Then $f_R:=f\circ\pi_R$ is in
$C^k_r(T\mathcal{D}^{s+k}_{\mu,D})$.
\end{theorem}

\begin{theorem} $(\pi_R$ is a Poisson map$)$\label{pi_R}
Let $k\geq 1$ and $r>\frac{1}{2}\operatorname{dim}M+2$ such that $s+k\geq
r$.
Then :
\[\{f\circ\pi_R,g\circ\pi_R\}^1(u_\eta)
=\left(\{f,g\}^1_+\circ\pi_R\right)
(u_\eta),\quad \forall\,f,g\in
C^k_r(\mathcal{V}^s_{\mu,D}),\quad u_\eta\in
T\mathcal{D}^{s+k}_{\mu,D}.\]
\end{theorem}
\textit{Proof} : Let $u_\eta\in T\mathcal{D}^{s+k}_{\mu,D}$
and $u:=\pi_R (u_\eta)$. The proof is a direct computation
using Lemmas \ref{1}, \ref{2}, \ref {vertical} and
\ref{horizontal}. Indeed, formula \eqref{poisson bracket
one}:
\[
\{f\circ\pi_R,g\circ\pi_R\}^1(u_\eta)
=\mathcal{G}^1(\eta)\left(\frac{\delta
f_R}{\delta\eta}(u_\eta),\frac
{\delta g_R}{\delta u}(u_\eta)\right)
-\mathcal{G}^1(\eta)\left(\frac{\delta
f_R}{\delta u}(u_\eta),\frac{\delta g_R}{\delta
\eta}(u_\eta)\right).
\]
So it suffices to compute the first term:
\begin{align*}
&\mathcal{G}^1(\eta)\left(\frac{\delta
f_R}{\delta\eta}(u_\eta),\frac {\delta g_R}{\delta
u}(u_\eta)\right)\\
&\qquad=\frac{1}{2}\mathcal{G}^1(\eta)
\Big{(}TR_\eta\Big{[}B^\alpha (u,\delta f(u))
- B^\alpha(\delta f(u),u)
+ \mathcal{P}_e\big{(}\mathcal{D}^\alpha (\delta f(u),u)
- \mathcal{D}^\alpha(u,\delta f(u))\big{)}\Big{]}, \\
&\qquad \qquad \qquad \qquad \qquad \qquad
 TR_\eta(\delta g(u))\Big{)}\\
&\qquad=\frac{1}{2}\left\langle B^\alpha (u,\delta
f(u))-B^\alpha(\delta f(u),u)+\mathcal{D}^\alpha (\delta
f(u),u)-\mathcal{D}^\alpha(u,\delta f(u)),\delta
g(u)\right\rangle_1\\
&\qquad=\frac{1}{2}\Big{(}\langle
B^\alpha(u,\delta f(u)),\delta g(u)\rangle_1-
\langle B^\alpha(\delta f(u),u),\delta g(u)\rangle_1\\
&\qquad\qquad\qquad
+\langle \mathcal{D}^\alpha(\delta f(u),u)+\nabla_{\delta
f
(u)}u,\delta g(u)\rangle_1-\langle\nabla_{\delta f(u)}u,\delta g(u)\rangle_1\\
&\qquad\qquad\qquad-\langle \mathcal{D}^\alpha(u,\delta
f(u))+\nabla_u\delta
f
(u),\delta g(u)\rangle_1+\langle\nabla_u\delta f(u),\delta g(u)\rangle_1\Big{)}
\\
&\qquad=\frac{1}{2}\Big{(}\langle u,\nabla_{\delta g(u)}\delta
f(u)\rangle_1-
\langle\delta f(u),\nabla_{\delta g(u)}u\rangle_1-\langle u,\nabla_{\delta f
(u)}\delta g(u)\rangle_1\\
&\qquad\qquad\qquad-\langle\nabla_{\delta f(u)}u,\delta g(u)
\rangle_1+\langle\delta f(u),\nabla_u\delta g(u)
\rangle_1+\langle\nabla_u\delta f(u),
\delta g(u)\rangle_1\Big{)},
\end{align*}
where we have used Lemmas \ref{1} and \ref{2} in the last
equality. Finally, after cancellation of several terms we
obtain:
\begin{align*}
&\mathcal{G}^1(\eta)\left(\frac{\delta
f_R}{\delta\eta}(u_\eta),\frac {\delta g_R}{\delta
u}(u_\eta)\right)-\mathcal{G}^1(\eta)\left(\frac{\delta
f_R}{\delta u}(u_\eta),\frac{\delta g_R}{\delta
\eta}(u_\eta)\right)\\ &\qquad\qquad=\langle
u,\nabla_{\delta g(u)}\delta f(u)\rangle_1-\langle
u,\nabla_{\delta f(u)}\delta
g(u)\rangle_1=\left(\{f,g\}^1_+\circ\pi_R\right)
(u_\eta).\,\,\blacksquare
\end{align*}

\begin{theorem} $(F_t$ is a Poisson map$)$\label{F_t}
Let $F_t$ be the flow of $\mathcal{S}^1$,
$t_1,t_2>\frac{1}{2}dim M+1$  such that $t_1\geq t_2$. Then
for all $G,H\in C^k_{t_2}(T\mathcal{D}^{t_1}_ {\mu,D})$ we
have:\\
$(1)$ $G\circ F_t, H\circ F_t\in C^k_{t_2}(T\mathcal{D}^{t_1}_{\mu,D})$\\
$(2)$ $\{G\circ F_t,H\circ F_t\}^1=\{G,F\}^1\circ F_t$ on $T\mathcal{D}^{t_1}_
{\mu,D}$.
\end{theorem}
\textit{Proof} : This is done as in Proposition 5.12 of
\cite {VaMa2004}. First of all we recall some general facts about weak
Riemannian Banach-manifolds. Let
$(Q,\langle\!\langle\,,\rangle\!\rangle)$ be a weak
Riemannian Banach-manifold with smooth geodesic spray. We define:
\[\mathcal{K}^\infty(TQ):=\left\{F\in C^\infty(TQ)  \Big|\exists\,\frac{\delta
F}{\delta \eta},\frac{\delta F}{\delta u}\in C^\infty(TQ,TQ)\right\}.
\]
Here, $\partial F/ \partial \eta$ and $\partial F/ \partial
u $ are the partial derivatives and $\delta F/ \delta
\eta$ and $\delta F / \delta u $ denote the horizontal and
vertical functional derivatives relative to the given
weak Riemannian metric on $Q $ of $F \in C^{\infty}(TQ)$
as defined at the beginning of section \ref{section: l p
structure}.

 Let $F_t$ be the geodesic flow and $G\in
\mathcal{K}^\infty(TQ)$. Then
$G\circ F_t\in\mathcal{K}^\infty(TQ)$ and $F_t$ is a Poisson map:
\begin{equation}
\label{Flow_is_Poisson}
\{G\circ F_t,H\circ F_t\}=\{G,H\}\circ F_t,\;
\forall\,G,H\in\mathcal{K}^\infty (TQ),
\end{equation}
where $\{\,,\}$ is the  Poisson bracket on
$\mathcal{K}^\infty(TQ)$ induced by the weak Riemannian
metric and the weak sympectic form on $T ^\ast Q$ (see
\eqref{poisson bracket one}).

We will use the following formula for $G\in\mathcal{K}^\infty(TQ)$:
\begin{equation}
\label{horizontal-vertical}
dG(u_\eta)(X_{u_\eta})=\frac{\partial G}{\partial \eta}(u_\eta)(T\pi_Q(X_
{u_\eta}))+\frac{\partial G}{\partial
u}(u_\eta)(K(X_{u_\eta})),
\end{equation}
where $\eta \in Q $, $u_\eta \in T_\eta Q $, $X_{u_\eta}
\in T_{u_\eta}(TQ)$, $\pi_Q: T Q \rightarrow Q $ is the
tangent bundle projection, and $K $ is the connector of
the given weak Riemannian metric on $Q$.
\medskip

With these general preparations, let $Q =
\mathcal{D}^{t_2}_{\mu,D}$ be endowed with the weak
Riemannian metric $\mathcal{G}^1$.
\medskip

\noindent (1) Let $G\in
C^k_{t_2}(T\mathcal{D}^{t_1}_{\mu,D})$, and
$u_\eta\in T\mathcal {D}^{t_1}_{\mu,D}$. So we have:
\[
\frac{\delta G}{\delta \eta}(F_t(u_\eta)),\frac{\delta G}{\delta u}(F_t
(u_\eta))\in T\mathcal{D}^{t_2}_{\mu,D}.
\]
Let $\widetilde{G}\in\mathcal{K}^\infty(T\mathcal{D}^{t_2}_{\mu,D})$ be such
that:
\[
\frac{\delta G}{\delta \eta}(F_t(u_\eta))=\frac{\delta \widetilde{G}}{\delta
\eta}(F_t(u_\eta))\;\text{ and }\;\frac{\delta G}{\delta u}(F_t(u_\eta))=\frac
{\delta \widetilde{G}}{\delta u}(F_t(u_\eta)).
\]
This is possible since  $\mathcal{D}^{t_2}_{\mu,D}$, and
hence $T\mathcal{D}^{t_2}_{\mu,D}$, are Hilbert manifolds
so they admit bump fuctions. Using
\eqref{horizontal-vertical} we find
\begin{align*}
&\frac{\partial (G\circ F_t)}{\partial
\eta}(u_\eta)=dG(F_t(u_\eta))\left(\frac
{\partial F_t}{\partial \eta}(u_\eta)\right)\\
&\quad=\mathcal{G}^1(u_\eta)\left(\frac{\delta G}{\delta \eta}(F_t
(u_\eta)),T\pi_{\mathcal{D}^s_{\mu,D}}\left(\frac{\partial F_t}{\partial \eta}
(u_\eta)\right)\right)+\mathcal{G}^1(u_\eta)\left(\frac{\delta G}{\delta u}(F_t
(u_\eta)),K^1\left(\frac{\partial F_t}{\partial \eta}(u_\eta)\right)\right)
\end{align*}
and so we obtain
\[
\frac{\partial (G\circ F_t)}{\partial \eta}(u_\eta)=\frac{\partial (\widetilde
{G}\circ F_t)}{\partial \eta}(u_\eta).
\]
 Since $\widetilde{G}\in\mathcal{K}^\infty(T\mathcal{D}^{t_2}_{\mu,D})$, we
obtain the existence of
\[
\frac{\delta (G\circ F_t)}{\delta \eta}(u_\eta)=\frac{\delta (\widetilde{G}
\circ F_t)}{\delta \eta}(u_\eta)\in T\mathcal{D}^{t_2}_{\mu,D}
\]
and the same is true for the vertical partial covariant
derivative. Doing this for  all $u_\eta\in
T\mathcal{D}^{t_1}_{\mu,D}$ we obtain that $G\circ F_t$ is
in
$C^k_{t_2}(T\mathcal{D}^{t_1}_{\mu,D})$.
\medskip

\noindent (2) Let $u_\eta$ be in
$T\mathcal{D}^{t_1}_{\mu,D}$. By part one, $\{G\circ
F_t,H\circ F_t\}^1(u_\eta)$ is well-defined and only
depends on
\[\frac{\delta G}{\delta \eta}(F_t(u_\eta)),\frac{\delta G}{\delta u}(F_t
(u_\eta)),\frac{\delta H}{\delta \eta}(F_t(u_\eta)),\frac{\delta H}{\delta u}
(F_t(u_\eta)).\]
Choosing $\widetilde{G}$ and $\widetilde{H}$ as in part one, and using \eqref
{Flow_is_Poisson} we obtain the desired formula.$\,\,\blacksquare$

\begin{theorem} $(\tilde{F}_t$ is a Poisson map$)$\label{tildeF_t}
Let $\tilde{F}_t=\pi_R\circ F_t$ be the flow of LAE-$\alpha$ equation. Then
we have
\[
\{f\circ\tilde{F}_t,g\circ\tilde{F}_t\}^1_+(u)
=\left(\{f,g\}^1_+\circ\tilde{F}_t\right)
(u),\qquad \forall f,g\in C^k_r(\mathcal{V}^s_{\mu,D}),
\qquad u\in\mathcal{V}^{s+2k}_{\mu,D},
\]
where $k\geq 1$ and $r>\frac{1}{2}dim M+2$ such that $s+k\geq r$ (for
example
$k=1$).
\end{theorem}
\textit{Proof} : Let $f\in C^k_r(\mathcal{V}^s_{\mu,D})$.
We have $f\circ\pi_R\in C^k_r(T\mathcal{D}^{s+k}_{\mu,D})$
by Theorem  \ref{f_R}. Therefore, by part (1) of Theorem
\ref{F_t} we get $f\circ\pi_R\circ F_t\in
C^k_r(T\mathcal{D}^{s+k}_{\mu,D})$ and hence
$f\circ\tilde{F}_t=f\circ\pi_R\circ
F_t|_{\mathcal{V}^{s+k}_{\mu,D}}
\in C^k_r(\mathcal{V}^{s+k}_{\mu,D})$. Since $\pi_R(u) =
u$, we have
\begin{align*}
\{f\circ\tilde{F}_t,g\circ\tilde{F}_t\}^1_+(u)&=\{f\circ\tilde{F}
_t\circ\pi_R,g\circ\tilde{F}_t\circ\pi_R\}^1(u)\text{ by Theorem \ref{pi_R}}\\
&=\{f\circ\pi_R\circ F_t,g\circ\pi_R\circ F_t\}^1(u)\text{ by Proposition
\ref{commutative diagram}}\\
&=\{f\circ\pi_R,g\circ\pi_R\}^1\left(F_t(u)\right)\text{ by
Theorem \ref{F_t}}\\
&=(\{f,g\}^1_+\circ\pi_R)(F_t(u))\text{ by Theorem
\ref{pi_R}}\\
&=(\{f,g\}^1_+\circ\tilde{F}_t)(u) \text{ by Proposition
\ref{commutative diagram}}.
\end{align*}
Note that for the first equality we need $u\in\mathcal{V}^{(s+k)+k}_
{\mu,D}$ by Theorem \ref{pi_R}.\,\,$\blacksquare$\\

The last Theorem gives the Poisson formulation of the
LAE-$\alpha$ equation. We  recall that an integral curve
$u(t)$ of the LAE-$\alpha$ (or the Euler)  equation is
$C^1$ as a map in $\mathcal{V}^{s-1}_{\mu,D}$, but it is
believed  to be continuous but not differentiable as a map
in $\mathcal{V}^s_{\mu,D}$.

\begin{theorem}
\label{l p formulation}
Let $u(t)\subset\mathcal{V}^s_{\mu,D}$ be a
curve such that $u\in  C^0(I,\mathcal{V}^s_{\mu,D})\cap
C^1(I,\mathcal{V}^{s-1}_{\mu,D})
$. Then
\[\frac{d}{dt}f(u(t))=\{f,h\}^1_+(u(t)), \forall\,f\in C^1_s(\mathcal{V}^{s-1}_
{\mu,D})\Longleftrightarrow u(t) \text{ is a solution of LAE-$\alpha$ equation}
\]
where $h(u):=\frac{1}{2}\langle u,u\rangle_1$ is the reduced
Hamiltonian.
\end{theorem}
\textit{Proof} : We remark that $h\in C^1_s(\mathcal{V}^s_{\mu,D})$ with
$\delta h(u)=u$. We find:
\begin{align*}
\frac{d}{dt}f(u(t))&=Df(u(t))(\partial_tu(t))\\
&=\langle\delta f(u(t)),\partial_tu(t)\rangle_1
\end{align*}
and, by Lemma \ref{1},
\begin{align*}
\{f,h&\}^1_+(u(t))=\langle u(t),\nabla_{u(t)}
\delta f(u(t))\rangle_1-\langle u(t), \nabla_{\delta
f(u(t))}u(t)\rangle_1\\
&=-\langle\nabla_{u(t)}u(t)+\mathcal{D}^\alpha(u(t),u(t)),\delta
f(u(t))
\rangle_1-\langle(1-\alpha^2\Delta_r)u(t),
\nabla_{\delta f(u(t))}u(t)\rangle_0.
\end{align*}
Using the remarkable fact that $\nabla u^t\cdot\Delta_ru$ is in $\mathfrak{X}
^{s-2}$ (Lemma \ref{F-alpha}), and the identity $\nabla u^t\cdot u=\operatorname
{grad}(g(u,u))$, we obtain for the second term:
\begin{align*}
\langle(1-&\alpha^2\Delta_r)u(t),\nabla_{\delta f(u(t))}u(t)\rangle_0\\
&=\langle \nabla u(t)^t\cdot (1-\alpha^2\Delta_r)u(t),\delta
f(u(t))\rangle_0\\
&=\langle (1-\alpha^2\mathcal{L})^{-1}\nabla u(t)^t\cdot
(1-\alpha^2\Delta_r)u
(t),\delta f(u(t))\rangle_1\\
&=\langle (1-\alpha^2\mathcal{L})^{-1}\operatorname{grad}[g(u(t),u(t))]-(1-
\alpha^2\mathcal{L})^{-1}\alpha^2\nabla u(t)^t\cdot\Delta_ru(t),\delta f(u(t))
\rangle_1\\
&=-\langle (1-\alpha^2\mathcal{L})^{-1}\alpha^2\nabla u(t)^t\cdot\Delta_ru
(t),\delta f(u(t))\rangle_1\quad\text{by the Stokes
decomposition. }
\end{align*}
So we obtain by Lemma \ref{F-alpha}:
\begin{align*}
\{f,&h\}^1_+(u(t))\\
&=\langle-\nabla_{u(t)}u(t)-\mathcal{D}^\alpha(u(t),u(t))+(1-\alpha^2\mathcal
{L})^{-1}\alpha^2\nabla u(t)^t\cdot\Delta_ru(t),\delta f(u(t))\rangle_1\\
&=\langle-\nabla_{u(t)}u(t)-\mathcal{F}^\alpha(u(t))-(1-\alpha^2\mathcal{L})^{-
1}\alpha^2\operatorname{grad}(F(u(t))),\delta f(u(t))\rangle_1\\
&=-\langle\mathcal{P}_e\left(\nabla_{u(t)}u(t)+\mathcal{F}^\alpha(u(t))
\right),\delta f(u(t))\rangle_1.
\end{align*}
Thus $\displaystyle\frac{d}{dt}f(u(t))=\{f,h\}^1_+(u(t)), \forall\,f\in C^1_s
(\mathcal{V}^{s-1}_
{\mu,D})$ is equivalent to:
\[\partial_tu(t)+\mathcal{P}_e\left(\nabla_{u(t)}u(t)+\mathcal{F}^\alpha(u(t))
\right)=0\]
which is LAE-$\alpha.\,\,\blacksquare$

\section{The case of free-slip and mixed boundary conditions}

In this section we shall generalize all our results to the case of
free-slip and mixed boundary conditions. Note that setting
$\Gamma_1=\varnothing$ in the mixed case, gives the free-slip case. The
fundamental difference between these boundary conditions and the no-slip
case we studied before is the following. For all vector fields $u,v$ in
$\mathcal{V}^s_D$, the vector field $\nabla_uv$ lies in
$\mathcal{V}^{s-1}_D$. This is a fact we used several times in our
previous computations. Unfortunately, for vector fields $u,v$ in
$\mathcal{V}^s_{mix}$ this is not true since $\nabla_uv$ may not be in
$\mathcal{V}^{s-1}_{mix}$. In this case we will use that
$\nabla_uv-\nabla_vu=[u,v]$ is in  $\mathcal{V}^{s-1}_{mix}$. As a first
consequence, the useful identity \eqref{lemma3} for the no-slip case
\[
(1-\alpha^2\mathcal{L})^{-1}\nabla_u[(1-\alpha^2\Delta_r)v]
=\nabla_uv+\mathcal{D}^{\alpha}(u,v),
\]
where $u$ is in $\mathcal{V}^s_{\mu,D},
s>\frac{1}{2}\operatorname{dim}M+1$, $v$ is in $\mathcal{V}^r_{\mu,D}$,
and $r>\frac{1}{2}\operatorname{dim}M+3$, is replaced by
\begin{equation}
\label{newlemma3}
(1-\alpha^2\mathcal{L})^{-1}\nabla_u[(1-\alpha^2\Delta_r)v]=(1-\alpha^2\mathcal{L})^{-1}(1-\alpha^2\mathcal{L})\nabla_uv+\mathcal{D}^{\alpha}(u,v)
\end{equation}
if $u$ is in $\mathcal{V}^s_{\mu,mix},
s>\frac{1}{2}\operatorname{dim}M+1$,
$v$ is in $\in\mathcal{V}^r_{\mu,mix}$, and
$r>\frac{1}{2}\operatorname{dim}M+3$.

Recall that for $r\geq 1$, $(1-\alpha^2\mathcal{L})$ denotes the continuous linear map $(1-\alpha^2(\Delta + 2\operatorname{Ric}+\operatorname{grad}\operatorname{div})) : \mathfrak{X}^r\longrightarrow \mathfrak{X}^{r-2}$ acting on all $H^r$ vector fields, and $(1-\alpha^2\mathcal{L})^{-1} : \mathfrak{X}^{r-2}\longrightarrow\mathcal{V}^r_{mix}$ denotes the inverse of the isomorphism $(1-\alpha^2\mathcal{L})_{|\mathcal{V}^r_{mix}}$.
Formula \eqref{newlemma3} induces some changes in Lemmas \ref{1} and \ref{2} which must be replaced by the following.

\begin{lemma}
\label{new1}
Let $s>\frac{1}{2}\operatorname{dim}M+1$. Let $u,v\in\mathcal{V}^s_{\mu,mix}$ and $w\in\mathcal{V}^s_{mix}$. Then:
\[\langle (1-\alpha^2\Delta_r)v,\nabla_uw\rangle_0=-\langle(1-\alpha^2\mathcal{L})^{-1}(1-\alpha^2\mathcal{L})\nabla_uv+\mathcal{D}^\alpha
(u,v),w\rangle_1\]
where $\mathcal{D}^\alpha : \mathcal{V}^s_{\mu,mix}\times\mathcal{V}^s_{\mu,mix}
\longrightarrow\mathcal{V}^s_{mix}$ is the bilinear continuous map given by
\begin{align*}
\mathcal{D}^\alpha(u,v)&:=\alpha^2(1-\alpha^2\mathcal{L})^{-1}\Big{(}
\operatorname{Div}(\nabla v\cdot\nabla u^t+\nabla v\cdot\nabla u)\\
&\qquad
\qquad+\operatorname{Tr}\big{(}\nabla_{\cdot}(\operatorname{R}(\cdot,u)
v)+\operatorname{R}(\cdot,u)\nabla_\cdot v\big{)}\\
&\qquad \qquad+\operatorname{grad}\big{(}\operatorname{Tr}(\nabla
u\cdot\nabla
v)+\operatorname{Ricci}(u,v)\big{)}-(\nabla_u\operatorname{Ric})(v)\Big{)}
\end{align*}

\end{lemma}
\textit{Proof} : Using the first part of Lemma \ref{1} and formula \eqref{newlemma3} we obtain for $u\in\mathcal{V}^s_{\mu,mix},w\in\mathcal{V}^s_{mix}$ and $v\in\mathcal{V}^r_{\mu,mix},r>\frac{1}{2}\operatorname{dim}M+3$:
\begin{eqnarray*}
\langle(1-\alpha^2\Delta_r)v,\nabla_uw\rangle_0&=&-\langle\nabla_u[(1-\alpha^2\Delta_r)v],w\rangle_0\\
&=&-\langle (1-\alpha^2\mathcal{L})^{-1}\nabla_u[(1-\alpha^2\Delta_r)
v],w\rangle_1\\
&=&-\langle(1-\alpha^2\mathcal{L})^{-1}(1-\alpha^2\mathcal{L})\nabla_uv+\mathcal{D}^{\alpha}(u,v),
w\rangle_1.
\end{eqnarray*}
Using the fact that
$\mathcal{V}^r_{\mu,mix},r>\frac{1}{2}\operatorname{dim}
M+3$ is dense in $\mathcal{V}^s_{\mu,mix}$, and the fact that
$\langle\,,\rangle_1, \nabla$, and $\mathcal{D}^\alpha$ are continuous on
$\mathcal{V}^s_{\mu,D}$, and $(1-\alpha^2\mathcal{L})^{-1}(1-\alpha^2\mathcal{L})$ is continuous on
$\mathcal{V}^{s-1}_{\mu,D}$ we obtain the desired result.$\,\, \blacksquare$

\begin{lemma} \label{new2}
Let $s>\frac{1}{2}\operatorname{dim} M+1$.
Let $B^\alpha : \mathcal{V}^{s+1}_{\mu,mix}\times\mathfrak{X}
^s\longrightarrow\mathcal{V}^{s+1}_{\mu,mix}$ the continuous bilinear map given
by
\[B^\alpha(v,w):=\mathcal{P}_e(1-\alpha^2\mathcal{L})^{-1}(\nabla w^t\cdot(1-
\alpha^2\Delta_r)v).\]
Then we have
\[\langle (1-\alpha^2\Delta_r)v,\nabla_uw\rangle_0=\langle B^\alpha(v,w),u\rangle_1\]
for all $u\in\mathcal{V}^r_{\mu,mix},r>\frac{1}{2}\operatorname{dim}M$, and
for all $v\in\mathcal{V}^{s+1}_{\mu,mix},$ and $w\in\mathfrak{X}^s$.
\end{lemma}
\textit{Proof} : The proof is similar to that of Lemma \ref{1}. Note that $\langle (1-\alpha^2\Delta_r)v,\nabla_uw\rangle_0$ does not equal $\langle v,\nabla_uw\rangle_1$ since $\nabla_uw$ does not belong to $\mathcal{V}^r_{mix}$.$\,\, \blacksquare$
\medskip

In order to carry out the Lie-Poisson reduction procedure for the mixed boundary conditions, we have to establish the existence and the smoothness of the geodesic spray of the weak Riemannian manifold $(\mathcal{D}^s_{\mu,mix},\mathcal{G}^1)$. So we will need a reformulation of LAE-$\alpha$ similar to \eqref{LAE} in the case of mixed boundary conditions. This reformulation is given by the following proposition where we use the Euler-Poincar\'{e} reduction theorem.

\begin{proposition}
Let $\eta(t)$ be a curve in $\mathcal{D}^s_{\mu,mix}$, and let $u(t):=TR_{\eta(t)^{-1}}(\dot\eta(t))=\dot\eta(t)\circ\eta(t)^{-1}
\in\mathcal{V}^s_{\mu,mix}$.
Then the following properties are equivalent :\\
$(1)$ $\eta(t)$ is a geodesic of $(\mathcal{D}^s_{\mu,mix},\mathcal{G}^1)$\\
$(2)$ $u(t)$ is a solution of :
\begin{equation}
\label{LAEnew}
\partial_tu(t)+\mathcal{P}_e\left((1-\alpha^2\mathcal{L})^{-1}(1-\alpha^2\mathcal{L})\nabla_{u(t)}u(t)
+\mathcal{F}^{\alpha}(u(t))\right)=0
\end{equation}
\end{proposition}
\textit{Proof} : By the the Euler-Poincar\'{e} reduction theorem, $\eta(t)$ is a geodesic of $(\mathcal{D}^s_{\mu,mix},\mathcal{G}^1)$ if and only if $u(t):=\dot\eta(t)\circ\eta(t)^{-1}$ is an extremum of the reduced action
\[
s(u)=\frac{1}{2}\int_a^b\langle u(t),u(t)\rangle_1dt
\]
for variations of the form
\[
\delta u(t)=\partial_tw(t)+[u(t),w(t)]
\]
where $w(t):=\delta \eta(t)\circ\eta^{-1}(t)$ vanishes at the endpoints. Integrating by parts, using the fact that $[u(t),w(t)]$
is in $\mathcal{V}^{s-1}_{\mu,mix}$ and with Lemma \ref{new2} we find:
\begin{align*}
&Ds(u)(\delta u)=\int_a^b\langle u(t),\delta u(t)\rangle_1dt\\
&=\int_a^b\langle u(t),\partial_tw(t)\rangle_1dt+\int_a^b\langle u(t),[u(t),w(t)]\rangle_1dt\\
&=-\int_a^b\langle \partial_tu(t),w(t)\rangle_1dt+\int_a^b\langle (1-\alpha^2\Delta_r)u(t),[u(t),w(t)]\rangle_0dt\\
&=-\int_a^b\langle \partial_tu(t),w(t)\rangle_1dt+\int_a^b\langle (1-\alpha^2\Delta_r)u(t),\nabla_{u(t)}w(t)\rangle_0dt\\
&\qquad-\int_a^b\langle (1-\alpha^2\Delta_r)u(t),\nabla_{w(t)}u(t)\rangle_0dt\\
&=-\int_a^b\langle \partial_tu(t),w(t)\rangle_1dt-\int_a^b\langle (1-\alpha^2\mathcal{L})^{-1}(1-\alpha^2\mathcal{L})\nabla_{u(t)}u(t)+\mathcal{D}^{\alpha}(u(t),u(t)),w(t)\rangle_1dt\\
&\qquad-\int_a^b\langle \nabla u(t)^t\cdot
(1-\alpha^2\Delta_r)u(t),w(t)\rangle_0dt.
\end{align*}

With Lemma \ref{1} (1), we have $\langle \nabla u(t)^t\cdot u(t),w(t)\rangle_0=0$, thus the last term equals
\[\alpha^2\int_a^b\langle \nabla u(t)^t\cdot \Delta_r u(t),w(t)\rangle_0dt \]
and we obtain:
\begin{align*}
Ds(u)(\delta u) & =
-\int_a^b\left\langle \partial_tu(t)
+(1-\alpha^2\mathcal{L})^{-1}(1-\alpha^2\mathcal{L})\nabla_{u(t)}u(t)
\right.\\
& \qquad +\left. \mathcal{D}^{\alpha}(u(t),u(t))
-\alpha^2(1-\alpha^2\mathcal{L})^{-1}\nabla
u(t)^t\cdot \Delta_r u(t),w(t)\right\rangle_1dt.
\end{align*}
So by the Stokes decomposition theorem, $Ds(u)(\delta u)=0$ for all
$\delta u$, is equivalent to
\begin{align*}
\mathcal{P}_e(\partial_tu(t)
&+(1-\alpha^2\mathcal{L})^{-1}(1-\alpha^2\mathcal{L})\nabla_{u(t)}u(t)\\
&+\mathcal{D}^{\alpha}(u(t),u(t))-\alpha^2(1-\alpha^2\mathcal{L})^{-1}\nabla
u(t)^t\cdot \Delta_r u(t))
=0
\end{align*}
and, with Lemma \ref{Falpha} (which remains valid on
$\mathcal{V}^s_{\mu,mix}$), this is equivalent to
\[
\mathcal{P}_e\big(\partial_tu(t)+(1-\alpha^2\mathcal{L})^{-1}(1-\alpha^2\mathcal{L})\nabla_{u(t)}u(t)+\mathcal{F}^{\alpha}(u(t))\big)=0.\,\, \blacksquare
\]

Let $\eta\in\mathcal{D}^s_{mix}$, $r\geq 0$, and $H^{r}_\eta:=\{u_\eta\in H^r(M,TM)|\pi\circ u=\eta\}$. We denote by $H^{r}_\eta\downarrow \mathcal{D}^s_{mix}$ the vector bundle over $\mathcal{D}^s_{mix}$, whose fiber at $\eta\in\mathcal{D}^s_{mix}$ is $H^{r}_\eta$.
The proof of Proposition 5 in \cite{Shkoller2000} shows that for $s>\frac{1}{2}\operatorname{dim}M+1$, the map
\[
\overline{(1-\alpha^2\mathcal{L})} : H^{s}_\eta\downarrow \mathcal{D}^s_{mix}\longrightarrow H^{s-2}_\eta\downarrow \mathcal{D}^s_{mix}
\]
defined by $\overline{(1-\alpha^2\mathcal{L})}(u_\eta):=[(1-\alpha^2\mathcal{L})(u_\eta\circ\eta^{-1})]\circ\eta$ is a $C^{\infty}$ bundle map. Furthermore,
\[
\overline{(1-\alpha^2\mathcal{L})} : T\mathcal{D}^s_{mix}\longrightarrow H^{s-2}_\eta\downarrow \mathcal{D}^s_{mix}
\]
is a bijection, whose inverse is denoted by
\[
\overline{(1-\alpha^2\mathcal{L})}^{-1} : H^{s-2}_\eta\downarrow \mathcal{D}^s_{mix}\longrightarrow T\mathcal{D}^s_{mix}
\]

With the same method and notations as in section 3, but using equation \eqref{LAEnew} instead of \eqref{LAE}, we obtain the following lemma

\begin{lemma}
\label{newspray} The geodesic spray of $(\mathcal{D}^s_{\mu,mix},\mathcal{G}^1)$ is given by:
\[
\mathcal{S}^1(u_\eta)=T\overline{\mathcal{P}}\left[T\left(\overline{(1-\alpha^2\mathcal{L})}^{-1}\circ \overline{(1-\alpha^2\mathcal{L})}\right)(S\circ u_\eta)-\operatorname{Ver}_{u_\eta}(\overline{\mathcal{F}}^{\alpha}(u_\eta))\right],
\]
where $S$ is the geodesic spray of $(M,g)$.

The connector $K^1 : TT\mathcal{D}^s_{\mu,mix}\longrightarrow T\mathcal{D}^s_{\mu,mix}$ of $(\mathcal{D}^s_{\mu,mix},\mathcal{G}^1)$ is given by:
\[
K^1(X_{u_\eta})=\overline{\mathcal{P}}\left(\overline{(1-\alpha^2\mathcal{L})}^{-1}\circ \overline{(1-\alpha^2\mathcal{L})}(K\circ X_{u_\eta})+\overline{\mathfrak{F}}^{\alpha}\big{(}\pi_{_{T\mathcal{D}^s_{\mu,mix}}}(X_{u_\eta}),T\pi_{_{\mathcal{D}^s_{\mu,mix}}}(X_{u_\eta})\big{)}\right),
\]
where $K :  TTM\longrightarrow TM$ is the connector of $(M,g)$.
\end{lemma}

Because of the existence of the geodesic spray $\mathcal{S}^1\in\mathfrak{X}^{C^{\infty}}(T\mathcal{D}^s_{\mu,mix})$ of the weak Riemannian manifold $(\mathcal{D}^s_{\mu,mix},\mathcal{G}^1)$, we can define the sets $C^k_r(T\mathcal{D}^t_{\mu,mix})$, the Poisson bracket $\{\,\,,\,\}^1$ on $C^k_r(T\mathcal{D}^t_{\mu,mix})$, the sets $C^k_{r,t}(\mathcal{V}^s_{\mu,mix})$ and $\mathcal{K}^k_{r,t}(\mathcal{V}^s_{\mu,mix})$, and the Poisson bracket $\{\,\,,\,\}^1_+$ on $C^k_{r,t}(\mathcal{V}^s_{\mu,mix})$ exactly in the same way we did in the case of no-slip boundary conditions.

As we shall see, all the properties of the Poisson bracket $\{\,\,,\,\}^1_+$ on
$C^k_{r,t}(\mathcal{V}^s_{\mu,mix})$ (Theorem \ref{properties1} and
\ref{properties2}) are still true in the mixed case but since the
Levi-Civita connection does not preserve the boundary conditions, the
computations in the proofs are more subtle.

\begin{theorem} \label{properties1new}
Let $s>\frac{1}{2}\operatorname{dim}M+1$ and $k\geq 1$. Then:
\[\{\,\,,\,\}^1_+:\mathcal{K}^k(\mathcal{V}^s_{\mu,mix})\times\mathcal{K}^k
(\mathcal{V}^s_{\mu,D})\longrightarrow\mathcal{K}^{k-1}_{s+1,s-1}(\mathcal{V}
^s_{\mu,mix})\]
and for all $u\in\mathcal{V}^{s+1}_{\mu,mix}$ we have
\begin{align*}
&\delta(\{f,g\}^1_+)(u)=\mathcal{P}_e(\nabla_{\delta g(u)}\delta
f(u)-\nabla_
{\delta f(u)}\delta g(u))\\
&\qquad+D\delta g(u)\left(\mathcal{P}_e\left((1-\alpha^2\mathcal{L})^{-1}(1-\alpha^2\mathcal{L})\nabla_{\delta
f(u)}u+\mathcal{D}
^\alpha(\delta f(u),u)\right)+B^\alpha(u,\delta f(u)\right)\\
&\qquad-D\delta f(u)\left(\mathcal{P}_e\left((1-\alpha^2\mathcal{L})^{-1}(1-\alpha^2\mathcal{L})\nabla_{\delta
g(u)}u+\mathcal{D}
^\alpha(\delta g(u),u)\right)+B^\alpha(u,\delta g(u)\right).
\end{align*}
\end{theorem}
\textit{Proof} : Let $h:=\{f,g\}^1_+$. We have to show that
$h\in\mathcal{K}^{k-1}_{s+1,s-1}(\mathcal{V}^s_{\mu,mix})$. As in Theorem
\ref{properties1} we obtain that $h\in C^k(\mathcal{V}^s_{\mu,mix})$, so
we can compute
$Dh(u)(v)$. Let $u,v\in\mathcal{V}^{s+1}_{\mu,mix}$. Using Lemmas
\ref{useful}, \ref{new1}, \ref{new2}, and \ref{3} (still valid in the
mixed case) we obtain:
\begin{align*}
&Dh(u)(v)\\
&=\langle v,\nabla_{\delta g(u)}\delta f(u)\rangle_1+\langle
u,\nabla_{D\delta g(u)(v)}\delta f(u)\rangle_1 +\langle u,\nabla_{\delta
g(u)}D\delta f(u)(v)\rangle_1\\
&\qquad-\langle v,\nabla_{\delta
f(u)}\delta g(u)\rangle_1-\langle u,\nabla_{D\delta f(u)(v)}\delta
g(u)\rangle_1-\langle u,\nabla_{\delta f(u)}D\delta g(u)(v)\rangle_1\\
&=\langle v,[\delta g(u),\delta f(u)]\rangle_1+\langle u,[D\delta
g(u)(v),\delta f(u)]\rangle_1 +\langle u,[\delta g(u),D\delta
f(u)(v)]\rangle_1\\
&=\langle v,[\delta g(u),\delta
f(u)]\rangle_1 + \langle (1-\alpha^2\Delta_r)u,[D\delta g(u)(v),\delta
f(u)]\rangle_0 \\
& \qquad + \langle (1-\alpha^2\Delta_r)u,[\delta
g(u),D\delta f(u)(v)]\rangle_0\\
&=\langle v,[\delta g(u),\delta
f(u)]\rangle_1+\langle (1-\alpha^2\Delta_r)u,\nabla_{D\delta
g(u)(v)}\delta f(u)\rangle_0 \\
& \qquad -\langle
(1-\alpha^2\Delta_r)u,\nabla_{\delta f(u)(v)}D\delta g(u)(v)\rangle_0
+\langle
(1-\alpha^2\Delta_r)u,\nabla_{\delta g(u)}D\delta
f(u)(v)\rangle_0 \\
&\qquad -\langle (1-\alpha^2\Delta_r)u,\nabla_{D\delta
f(u)(v)}\delta g(u)\rangle_0\\
&=\langle v,[\delta g(u),\delta
f(u)]\rangle_1+\langle B^\alpha(u,\delta f(u)),D\delta g(u)(v)\rangle_1\\
&\qquad+\langle
(1-\alpha^2\mathcal{L})^{-1}(1-\alpha^2\mathcal{L})\nabla_{\delta
f(u)}u+\mathcal{D}^\alpha(\delta f(u),u),D\delta g(u)(v)\rangle_1\\
&\qquad-\langle
(1-\alpha^2\mathcal{L})^{-1}(1-\alpha^2\mathcal{L})\nabla_{\delta
g(u)}u+\mathcal{D}^\alpha(\delta g(u),u),D\delta f(u)(v)\rangle_1\\
&\qquad-\langle B^\alpha(u,\delta g(u)),D\delta f(u)(v)\rangle_1\\
&=\langle v,[\delta g(u),\delta f(u)]\rangle_1+\langle D\delta
g(u)\left(B^\alpha(u,\delta f(u)\right),v\rangle_1\\ &\qquad+\langle
D\delta
g(u)\left(\mathcal{P}_e\left((1-\alpha^2\mathcal{L})^{-1}(1-\alpha^2\mathcal{L})\nabla_{\delta
f(u)}u+\mathcal{D}^\alpha(\delta f(u),u)\right)\right),v\rangle_1\\
&\qquad-\langle D\delta
f(u)\left(\mathcal{P}_e\left((1-\alpha^2\mathcal{L})^{-1}(1-\alpha^2\mathcal{L})\nabla_{\delta
g(u)}u+\mathcal{D}^\alpha(\delta g(u),u)\right)\right),v\rangle_1\\
&\qquad-\langle D\delta f(u)\left(B^\alpha(u,\delta
g(u)\right),v\rangle_1.
\end{align*}
Now the result follows as in Theorem \ref{properties1}.$\,\, \blacksquare$

\begin{theorem}
\label{properties2new}
Let $s,t>\frac{1}{2}\operatorname{dim} M+1,\,r\geq s,$ and $k\geq 1$.\\
$(1)$ $\{\,\,,\,\}^1_+$ is $\mathbb{R}$-bilinear and anti-symmetric on $C^k_
{r,t}(\mathcal{V}^s_
{\mu,mix})\times C^k_{r,t}(\mathcal{V}^s_{\mu,mix})$.\\
$(2)$ $\{\,\,,\,\}^1_+$ is a derivation in each factor:
\[\{fg,h\}_+^1=\{f,h\}_+^1g+f\{g,h\}_+^1,\forall\,f,g,h\in C^k_{r,t}(\mathcal
{V}^s_{\mu,mix}).\]
$(3)$ If $s>\frac{1}{2}\operatorname{dim} M+2$, $\{\,\,,\,\}^1_+$satisfies
the Jacobi identity:\\
For all $f,g,h\in\mathcal{K}^k(\mathcal{V}^s_{\mu,mix})$ and $u\in \mathcal{V}^
{s+1}_{\mu,mix}$ we have:
\[\{f,\{g,h\}^1_+\}^1_+(u)+\{g,\{h,f\}^1_+\}^1_+(u)+\{h,\{f,g\}^1_+\}^1_+(u)=0
\]
\end{theorem}
\textit{Proof}: $(1)$ This is obvious.\\
$(2)$ Let $f,g,h\in C^k_{r,t}(\mathcal{V}^s_{\mu,mix})$, and $u\in\mathcal{V}^r_{\mu,mix}$. Using Lemmas \ref{useful} and \ref{new1} we find:
\begin{align*}
\{fg,h\}_+^1(u)&=\langle u,[\delta (fg)(u),\delta h(u)]\rangle_1\\
&=\langle (1-\alpha^2\Delta_r)u,\nabla_{\delta (fg)(u)}\delta h(u)\rangle_0-\langle (1-\alpha^2\Delta_r)u,\nabla_{\delta h(u)}\delta (fg)(u)\rangle_0\\
&=\langle (1-\alpha^2\Delta_r)u,\nabla_{\delta f(u)}\delta h(u)\rangle_0g(u)+\langle (1-\alpha^2\Delta_r)u,\nabla_{\delta g(u)}\delta h(u)\rangle_0f(u)\\
&\qquad+\langle (1-\alpha^2\mathcal{L})^{-1}(1-\alpha^2\mathcal{L})\nabla_{\delta h(u)}u+\mathcal{D}^\alpha(\delta
h(u),u),\delta (fg)(u)\rangle_1\\
&=\langle (1-\alpha^2\Delta_r)u,\nabla_{\delta f(u)}\delta h(u)\rangle_0g(u)+\langle (1-\alpha^2\Delta_r)u,\nabla_{\delta g(u)}\delta h(u)\rangle_0f(u)\\
&\qquad+\langle (1-\alpha^2\mathcal{L})^{-1}(1-\alpha^2\mathcal{L})\nabla_{\delta h(u)}u+\mathcal{D}^\alpha(\delta
h(u),u),\delta f(u)\rangle_1g(u)\\
&\qquad+\langle (1-\alpha^2\mathcal{L})^{-1}(1-\alpha^2\mathcal{L})\nabla_{\delta h(u)}u+\mathcal{D}^\alpha(\delta
h(u),u),\delta g(u)\rangle_1f(u)\\
&=\langle (1-\alpha^2\Delta_r)u,\nabla_{\delta f(u)}\delta h(u)\rangle_0g(u)+\langle (1-\alpha^2\Delta_r)u,\nabla_{\delta g(u)}\delta h(u)\rangle_0f(u)\\
&\qquad-\langle (1-\alpha^2\Delta_r)u,\nabla_{\delta h(u)}\delta f(u)\rangle_0g(u)-\langle (1-\alpha^2\Delta_r)u,\nabla_{\delta h(u)}\delta g(u)\rangle_0f(u)\\
&=\langle (1-\alpha^2\Delta_r)u,[\delta f(u),\delta h(u)]\rangle_0g(u)+\langle (1-\alpha^2\Delta_r)u,[\delta g(u),\delta h(u)]\rangle_0f(u)\\
&=\langle u,[\delta f(u),\delta h(u)]\rangle_1g(u)+\langle u,[\delta g(u),\delta h(u)]\rangle_1f(u)\\
&=\{f,h\}_+^1(u)g(u)+f(u)\{g,h\}_+^1(u).
\end{align*}
$(3)$ Let $f,g,h\in\mathcal{K}^k(\mathcal{V}^s_{\mu,mix})$,
and $u\in\mathcal {V}^{s+1}_{\mu,D}$.
By Theorem \ref{properties1} we obtain
$\{g,h\}^1_+\in\mathcal{K}^{k-1}_
{s+1,s-1}(\mathcal{V}^s_{\mu,mix})\subset
C^k_{s+1,s-1}(\mathcal{V}^s_{\mu,mix})$.  Since
$s-1>\frac{1}{2}\operatorname{dim} M+1$ we can compute
the expression
$\{f,\{g,h\}^1_+\}
^1_+(u)$. We have:
\begin{align*}
&\{f,\{g,h\}^1_+\}^1_+(u)\\
&\qquad=\langle u,[\delta \{g,h\}^1_+(u),\delta f(u)]\rangle_1\\
&\qquad=\langle (1-\alpha^2\Delta_r)u,\nabla_{\delta \{g,h\}^1_+(u)}\delta
f(u)\rangle_1-\langle (1-\alpha^2\Delta_r)u,\nabla_{\delta f(u)}\delta \{g,h\}^1_+(u)\rangle_1
\end{align*}
So we can use Lemmas \ref{new1} and \ref{new2} and then the expression
for $\delta \{g,h\}^1_+(u)$ in Theorem \ref{properties1new}. Doing
exactly the same computation as in Theorem \ref{properties2} and using
analogous notations we find
\begin{align*}
&\{f,\{g,h\}^1_+\}^1_+(u)\\
&=\langle
\nabla_{[\delta h(u),\delta g(u)]}\delta f(u),(1-\alpha^2\Delta_r)u\rangle_0-
\langle\nabla_{\delta f(u)}[\delta h(u),\delta g(u)],(1-\alpha^2\Delta_r)u\rangle_0+D_{hgf}-D_{ghf}\\
&=\langle[[\delta h(u),\delta g(u)],\delta
f(u)],(1-\alpha^2\Delta_r)u\rangle_0+D_{hgf}-D_{ghf}\\
&=\langle[[\delta h(u),\delta g(u)],\delta
f(u)],u\rangle_1+D_{hgf}-D_{gfh}.
\end{align*}
Using the Jacobi identity for the Jacobi-Lie bracket of vector fields we
obtain the desired result.$\,\, \blacksquare$
\medskip

As in section 5, for $f\in C^k_r(\mathcal{V}^s_{\mu,mix})$, we shall denote $f_R:=f\circ\pi_R\in C^k(T\mathcal{D}^{s+k}_{\mu,mix})$. The proof of Lemma \ref{vertical} remains valid in the mixed case, so if $f\in C^k_r(\mathcal{V}^s_{\mu,mix})$, $k\geq 1$ and $r>\frac{1}{2}\operatorname{dim}M+1$ are such that $s+k\geq r$, then the vertical functional derivative of $f_R$ with respect to $\mathcal{G}^1$ exists and is given by:
\[
\frac{\delta f_R}{\delta u}(u_\eta)=TR_\eta(\delta f(\pi_R(u_\eta)))\in
T\mathcal{D}^r_{\mu,D}, \qquad \forall\,u_\eta\in
T\mathcal{D}^{s+k}_{\mu,mix}.
\]
Lemma \ref{horizontal} about the horizontal functional derivative remains valid in the mixed case but some computations in the proof should be adapted to this case. These computations are given below.

\begin{lemma}\label{horizontalnew}
Let $k\geq 1$ and $r>\frac{1}{2}\operatorname{dim}M+2$ such that $s+k\geq
r$.
Let $f\in C^k_r (\mathcal{V}^s_{\mu,mix})$.
Then the horizontal functional derivative of $f_R$ with respect to $\mathcal{G}
^1$ exists. It is given by:
\[\frac{\delta f_R}{\delta \eta}(u_\eta)=\frac{1}{2}TR_h\Big{[}B^\alpha
(u,\delta f(u))-B^\alpha(\delta f(u),u)+\mathcal{P}_e\big{(}\mathcal{D}^\alpha
(\delta f(u),u)-\mathcal{D}^\alpha(u,\delta f(u))\big{)}\Big{]}\]
for all $u_\eta\in T\mathcal{D}^{s+k}_{\mu,mix}$, where $u:=\pi_R(u_\eta)$
and $B ^\alpha$ was defined in Lemma \ref{new2}. So
we have:
\[
\frac{\delta f_R}{\delta \eta}(u_\eta)\in
T\mathcal{D}^r_{\mu,mix}, \qquad \forall\,u_\eta\in
T\mathcal{D}^{s+k}_{\mu,mix}.
\]
\end{lemma}
\textit{Proof} : As in Lemma \ref{horizontal}, we find for $u_\eta,v_\eta\in T\mathcal{D}^{s+k}_{\mu,D}$:
\begin{align*}
\left\langle\frac{\partial f_R}{\partial
\eta}(u_\eta), v_\eta\right\rangle
&=-Df(u_\eta\circ\eta^{-1})
(K^1(T(u_\eta\circ \eta^{-1})\circ(v_\eta\circ \eta^{-1})))\\
&=-Df(u)(K^1(Tu\circ v))\qquad \text{where
$u:=u_\eta\circ
\eta^{-1}$ and
$v:=v_\eta\circ \eta^{-1}$}.
\end{align*}
With the formula for the connector in Lemma \ref{newspray} we obtain the following identity:
\[
K^1(Tu\circ
v)=\mathcal{P}_e((1-\alpha^2\mathcal{L})^{-1}(1-\alpha^2\mathcal{L})\nabla_vu+\mathfrak{F}^\alpha(u,v)).
\]
So using formula \eqref{F-alpha} (still valid in the mixed case), and the notation $L^\alpha:=(1-\alpha^2\mathcal{L})^{-1}(1-\alpha^2\mathcal{L})$ we find:
\begin{align*}
\left\langle\frac{\partial f_R}{\partial
\eta}(u_\eta), v_\eta\right\rangle
&=-\langle\delta
f(u),L^\alpha(\nabla_vu)+\frac{1}{2}(\mathcal{D}^\alpha(u,v)
+\mathcal{D}^\alpha(v,u))\rangle_1\\
&\qquad+\frac{1}{2}\langle\delta f(u),(1-\alpha^2\mathcal{L})^{-1}
\alpha^2\big{(}\operatorname{grad}(G(u,v))\big{)}\rangle_1\\
&\qquad+\frac{1}{2}\langle\delta
f(u),(1-\alpha^2\mathcal{L})^{-1}
\alpha^2\big{(}\nabla u^t\cdot\Delta_r v+\nabla v^t\cdot\Delta_r u\big{)}
\rangle_1\\
&=-\langle\delta
f(u),L^\alpha(\nabla_vu)+\frac{1}{2}(\mathcal{D}^\alpha(u,v)
+\mathcal{D}^\alpha(v,u))\rangle_1\\
&\qquad+\frac{1}{2}\langle\delta
f(u),\alpha^2\big{(}\nabla u^t\cdot\Delta_r v+\nabla v^t\cdot\Delta_r u\big{)}
\rangle_0\\
&=-\langle\delta
f(u),\frac{1}{2}L^\alpha(\nabla_vu)+\frac{1}{2}\mathcal{D}^\alpha(u,v)
+\frac{1}{2}(\mathcal{D}^\alpha(v,u)+L^\alpha(\nabla_vu))\rangle_1\\
&\qquad-\frac{1}{2}\langle\delta f(u),\nabla u^t\cdot(1-\alpha^2\Delta_r) v\rangle_0+\frac{1}{2}\langle\delta f(u),\nabla u^t\cdot v\rangle_0\\
&\qquad-\frac{1}{2}\langle\delta f(u),\nabla v^t\cdot(1-\alpha^2\Delta_r) u\rangle_0+\frac{1}{2}\langle\delta f(u),\nabla v^t\cdot
u\rangle_0.
\end{align*}
Using Lemmas \ref{new1} and \ref{new2}, we obtain:
\begin{align*}
&2\left\langle\frac{\partial f_R}{\partial
\eta}(u_\eta), v_\eta\right\rangle\\
&=-\langle\delta f(u),L^\alpha(\nabla_vu)\rangle_1-\langle\delta
f(u),\mathcal{D}^{\alpha}(u,v)\rangle_1+\langle
(1-\alpha^2\Delta_r)u,\nabla_v\delta f(u)\rangle_0\\
&\qquad
-\langle\nabla_{\delta
f(u)}u,(1-\alpha^2\Delta_r)v\rangle_0-\langle\nabla_{\delta
f(u)}v,(1-\alpha^2\Delta_r)u\rangle_0\\
& =-\langle\delta
f(u),L^\alpha(\nabla_vu)\rangle_1-\langle\delta
f(u),L^\alpha(\nabla_uv)+\mathcal{D}^{\alpha}(u,v)\rangle_1+\langle\delta
f(u),L^\alpha(\nabla_uv)\rangle_1\\
&\qquad +\langle
B^\alpha(u,\delta f(u)),v\rangle_1-\langle\nabla_{\delta
f(u)}u,(1-\alpha^2\Delta_r)v\rangle_0+\langle L^\alpha(\nabla_{\delta
f(u)}u)+\mathcal{D}^{\alpha}(\delta f(u),u),v\rangle_1\\
&=-\langle\delta f(u),L^\alpha(\nabla_vu)\rangle_1+\langle
(1-\alpha^2\Delta_r)v,\nabla_u\delta f(u)\rangle_0+\langle\delta
f(u),L^\alpha(\nabla_uv)\rangle_1 \\
&\qquad +\langle B^\alpha(u,\delta f(u)),v\rangle_1
-\langle\nabla_{\delta
f(u)}u,(1-\alpha^2\Delta_r)v\rangle_0+\langle L^\alpha(\nabla_{\delta
f(u)}u),v\rangle_1 \\
&\qquad +\langle \mathcal{D}^{\alpha}(\delta
f(u),u),v\rangle_1\\
&=\langle\delta
f(u),L^\alpha(\nabla_uv-\nabla_vu)\rangle_1+\langle
(1-\alpha^2\Delta_r)v,\nabla_u\delta f(u)-\nabla_{\delta
f(u)}u\rangle_0\\
&\qquad +\langle L^\alpha(\nabla_{\delta
f(u)}u),v\rangle_1+\langle B^\alpha(u,\delta f(u)),v\rangle_1+\langle
\mathcal{D}^{\alpha}(\delta f(u),u),v\rangle_1.
\end{align*}
Since the Jacobi-Lie bracket of vector fields preserves the mixed boundary condition we have:
\begin{align*}
&\langle\delta f(u),L^\alpha(\nabla_uv-\nabla_vu)\rangle_1+\langle (1-\alpha^2\Delta_r)v,\nabla_u\delta f(u)-\nabla_{\delta f(u)}u\rangle_0\\
& \qquad = \langle\delta f(u),\nabla_uv-\nabla_vu\rangle_1+\langle
v,\nabla_u\delta f(u)-\nabla_{\delta f(u)}u\rangle_1\\
& \qquad =\langle
(1-\alpha^2\Delta_r)\delta f(u),\nabla_uv-\nabla_vu\rangle_0+\langle
v,L^\alpha(\nabla_u\delta f(u)-\nabla_{\delta f(u)}u)\rangle_1\\
& \qquad =\langle (1-\alpha^2\Delta_r)\delta
f(u),\nabla_uv\rangle_0-\langle (1-\alpha^2\Delta_r)\delta
f(u),\nabla_vu\rangle_0 \\
&\qquad \qquad +\langle v,L^\alpha(\nabla_u\delta f(u))-\langle
v,L^\alpha(\nabla_{\delta f(u)}u)\rangle_1\\
& \qquad =-\langle L^\alpha(\nabla_u\delta
f(u))+\mathcal{D}^\alpha(u,\delta f(u),v\rangle_1-\langle B^\alpha(\delta
f(u),u),v\rangle_1 \\
&\qquad \qquad +\langle v,L^\alpha(\nabla_u\delta f(u))-\langle
v,L^\alpha(\nabla_{\delta f(u)}u)\rangle_1\\
& \qquad =-\langle
\mathcal{D}^\alpha(u,\delta f(u),v\rangle_1-\langle B^\alpha(\delta
f(u),u),v\rangle_1-\langle v,L^\alpha(\nabla_{\delta f(u)}u)\rangle_1.
\end{align*}
So we obtain
\[
2\left\langle\frac{\partial f_R}{\partial\eta}(u_\eta), v_\eta \right \rangle=\langle B^\alpha
(u,\delta f(u))-B^\alpha(\delta f(u),u)+\mathcal{D}^\alpha
(\delta f(u),u)-\mathcal{D}^\alpha(u,\delta f(u)),v\rangle_1
\]
and the result follow.$\,\, \blacksquare$
\medskip

We conclude that Theorem \ref{f_R} remains valid in the mixed case.
For proving that $\pi_R$ is a Poisson map (in the sense of Theorem
\ref{pi_R}) it suffices to use Lemmas \ref{new1} and \ref{new2}  insteed
of Lemmas \ref{1} and \ref{2} in the proof of Theorem \ref{pi_R}. So
Theorems \ref{F_t}, \ref{tildeF_t}, and \ref{l p formulation}  are
also valid in this case.

{\footnotesize

\bibliographystyle{new}

}

\end{document}